\documentclass[12pt,twoside]{article}
\usepackage{graphicx}
\usepackage{amsmath,amssymb}

\title{On a variational method for solving a class of $p$-Kirchhoff type problem} %
\author{ S.H. Rasouli , K. Fallah \\
 Department of Mathematics, Faculty of Basic Science,
\\Babol University of Technology, Babol, Iran\\e-mail:
s.h.rasouli@nit.ac.ir}

\date{}

\begin{document}
\maketitle
\begin{center}
{\bf\large Abstract}\\
\end{center}

In this paper we are concerned with some $p$-Kirchhoff type problems
involving sign-changing weight functions. We prove the existence of
multiple positive solutions of the problem via the Nehari
manifold approach.\\\\
\hspace{-0.6 cm}Keywords: $p$-Kirchhoff type problem; Multiple
solutions; Sign-changing weight functions; Nehari manifold.
\section{Introduction}
\hspace{0.6 cm}The aim of this paper is to prove some existence
and multiplicity results of solutions to the following problem:
\begin{equation}
\left\{\begin{array}{ll}
-M\,\Big(\int_{\Omega }(|\nabla
u|^{p})\,dx\Big) \Delta_{p}u=\lambda f(x)|u|^{q-2}u+g(x)|u|^{r-2}u, & x\in \Omega,\\
u=0  ,  & x\in\partial \Omega,
\end{array}\right.
\end{equation}
where $\Omega$ is a smooth bounded domain in $\mathbb{R}^N$,
$1<q<p<r\leq\,p^{*}$($p^{*}=\frac{pN}{N-p}$ if
$N\geq\,2p-1;p^{*}=\infty$ if $N=1,p$), $M(s)=as^{p-1}+b$ and
$a,b,\lambda
>0$. The weight functions $f,g\in C(\bar{\Omega})$ satisfy the
following conditions:
\begin{itemize}
\item[(i)] $f^{+}=\max\{f,0\}\neq 0$,
\item[(ii)] $g^{+}=\max\{g,0\}\neq 0$.\\
\end{itemize}

Problem $(1)$ is a general version of a model presented by Kirchhoff
\cite{gk}. More precisely, Kirchhoff introduced a model
\begin{equation}\label{E12}
\rho\frac{\partial^2u}{\partial
t^2}-\Big(\frac{\rho_0}h+\frac{E}{2L}\int_0^L|\frac{\partial
u}{\partial x}|^2dx\Big)\frac{\partial^2u}{\partial x^2}=0,
\end{equation}
where $\rho, \rho_0, h, E, L$ are constants, which extends the
classical D'Alembert's wave equation by considering the effects of
the changes in the length of the strings during the vibrations. The
 problem
\begin{equation}\label{E13}
\begin{gathered}
-\Big(a+b\int_\Omega|\nabla u|^2dx\Big)\Delta
u=f(x,u)\quad \text{in }\Omega\\
u=0\quad \text{on }\partial\Omega
\end{gathered}
\end{equation}
received much attention, mainly after the article by Lions
\cite{jll}. Problems like \eqref{E13} are also introduced as models
for other physical phenomena as, for example, biological systems
where $u$ describes a process which depends on the average of itself
(for example, population density). See \cite{ACM} and its references
therein. For a more detailed reference on this subject we refer the
interested reader to \cite{AP,CCS,CF,CN,M,PZ}.\\

See \cite{cyc-yck-tfw} where the authors discussed the problem $(1)$
when $p=2.$ Here we focus on extending the study in
\cite{cyc-yck-tfw}. In fact This paper is motivated, in part, by the
mathematical difficulty posed by the degenerate quasilinear elliptic
operator compared to the Laplacian operator ($p=2$). This extension
is nontrivial and requires more careful analysis of the
nonlinearity. Our approach is based on the Nehari manifold, see \cite{gaa-shr,kjb-tfw1,kjb-tfw2,kjb-yz,pd-sip,shr-gaa,tfw1,tfw2}.\\

\section{Variational setting}

\hspace{0.6 cm}We define\\
$$
\|u\|_L^{z}\leq \,S_{z}^{-\frac{z}{p}}\,\|u\|_W^{1,p},\,\,\text{for
all}\,\, u\in W_{0}^{1,p}(\Omega)\setminus\{0\},
$$
where $S_z$ is a best sobolev constant for the embedding of
$W_{0}^{1,p}(\Omega)$ in $L^{z}(\Omega)\,\,\,\text {with}\,\,
1<z<p*$.\\
The energy functional corresponding to equation $(1)$, for $u \in
W_{0}^{1}(\Omega) $ is defined
by\\
$$
\mathcal{J}_{\lambda,M}(u)=\frac{1}{p}\,\hat{M}(\|u\|_{W^{1,p}}^{p})
-\,\frac{\lambda}{q}\,\int_{\Omega}f|u|^{q}\,dx-\,\frac{1}{r}\,\int_{\Omega}g|u|^r\,dx,
$$
where $\hat{M}(s)=\int_{0}^{t} M(t)\,dt$ and $M:\mathbb{R}\to
\mathbb{R}^{+}$ is any function that is differentiable everywhere
except at some finite points. It is well known the weak solutions of
equation $(1)$ are the critical points of the energy functional
$\mathcal{J}_{\lambda,M}$. By taking $M(s)=as^{p-1}+b$ and using
the Sobolev inequality, we can write\\
\begin{eqnarray*}
\mathcal{J}_{\lambda,M}(u)&=&\frac{1}{p}\,\hat{M}(\|u\|_{W^{1,p}}^{p})
-\,\frac{\lambda}{q}\,\int_{\Omega}f|u|^{q}\,dx-\,\frac{1}{r}\,\int_{\Omega}g|u|^r\,dx\\
&\geq&\,\frac{a}{p^{2}}\,\|u\|_{W^{1,p}}^{p^{2}}+\,\frac{b}{p}\,\|u\|_{W^{1,p}}^{p}
-\frac{\lambda
S_{q}^{\frac{-q}{p}}\,\|f^{+}\|_{\infty}}{q}\,\|u\|_{W^{1,p}}^{q}
-\,\frac{S_{r}^{\frac{-r}{p}}\,\|g^{+}\|_{\infty}}{r}\,\|u\|_{W^{1,p}}^{r}\\
&=&\big(\frac{a}{p^{2}}\|u\|_{W^{1,p}}^{p^{2}-r}-\frac{S_{r}^{\frac{-r}{p}}\|g^{+}\|_{\infty}}{r}\Big)\|u\|_{W^{1,p}}^{r}
+\,\big(\frac{b}{p}\,\|u\|_{W^{1,p}}^{p-q}-\,\frac{\lambda
S_{q}^{\frac{-q}{p}}\,\|f^{+}\|_{\infty}}{q}\Big)\,\|u\|_{W^{1,p}}^{q},
\end{eqnarray*}
for all $u\in W_{0}^{1,p}\,(\Omega)\,\setminus\{0\}$, and so
$\mathcal{J}_{\lambda,M}\,(u)$ is bounded below on
$W_{0}^{1,p}\,(\Omega)$ when $r<p^{2}$ and, when $r>p^{2}$,
$\mathcal{J}_{\lambda,M}$ is no longer bounded below on
$W_{0}^{1,p}\,(\Omega)$, because
$\lim_{t\to\infty}\,\mathcal{J}_{\lambda,M}\,(tu)=\,-\,\infty$,
so it is useful to consider the functional on the Nehari manifold\\
\begin{eqnarray*}
\mathcal{N}_{\lambda,M}=\{u\in W_{0}^{1,p}\,(\Omega) \backslash
\{0\}:\,\langle \mathcal{J}_{\lambda,M }'(u) ,u\rangle =0\},
\end{eqnarray*}
where $\langle,\rangle$ denote the usual duality. Thus,$u\in
\mathcal{N}_{\lambda,M}$ if and only if
$$
M(\|u\|_{{W}^{1,p}}^{p})\,\|u\|_{{W}^{1,p}}^{p}-\,
\lambda\,\int_{\Omega}f|u|^{q}\,dx-\,\int_{\Omega}g|u|^{r}\,dx=0.
$$
The Nehari manifold $\mathcal{N}_{\lambda,M}$ is closely linked to
the behavior of functions of the form $I_{u,M}:\,t\to\
\mathcal{J}_{\lambda,M}(tu)$ for $t>0$ that named fibering maps.
If $u\in W_{0}^{1,p}$, we have\\
\begin{eqnarray*}
I_{u,M}(t)=\frac{1}{p}\,\hat{M}(t^{p}\,\|u\|_{{W}^{1,p}}^{p})
-\,\frac{\lambda
t^{q}}{q}\,\int_{\Omega}f|u|^{q}\,dx-\,\frac{t^{r}}{r}\,\int_{\Omega}g|u|^{r}\,dx,
\end{eqnarray*}
\begin{eqnarray*}
\hspace*{2cm}I_{u,M}'(t)=\,t^{p-1}\,M(t^{p}\,\|u\|_{{W}^{1,p}}^{p})\,\|u\|_{{W}^{1,p}}^{p}
-\,\lambda\,t^{q-1}\,\int_{\Omega}f|u|^{q}\,dx-\,t^{r-1}\,\int_{\Omega}g|u|^{r}\,dx,
\end{eqnarray*}
\begin{eqnarray*}
\hspace*{2cm}I_{u,M}''(t)&=&(p-1)\,t^{p-2}\,M(t^{p}\,\|u\|_{{W}^{1,p}}^{p})\,\|u\|_{{W}^{1,p}}^{p}
+\,p\,t^{2p-2}\,M'(t^{p}\,\|u\|_{{W}^{1,p}}^{p})\,\|u\|_{{W}^{1,p}}^{2p}\\
&-&\,\lambda\,(q-1)\,t^{q-2}\,\int_{\Omega}f|u|^{q}\,dx-\,(r-1)\,t^{r-2}\,\int_{\Omega}g|u|^{r}\,dx.
\end{eqnarray*}
Clearly,\\
\begin{eqnarray*}
\hspace*{2cm}t\,I_{u,M}'(t)=\,M(\|tu\|_{{W}^{1,p}}^{p})\,\|tu\|_{{W}^{1,p}}^{p}-\,\lambda\,\int_{\Omega}f|tu|^{q}\,dx
-\,\int_{\Omega}g|tu|^{r}\,dx;
\end{eqnarray*}
and so, for $u\in W_{0}^{1,p}(\Omega)\backslash\{0\}$ and $t>0$,
$I_{u,M}'(t)=0$ if and only if $tu\in \mathcal{N}_{\lambda,M}$,
i.e., positive critical points of $I_{u,M}$ correspond to points on
the Nehari manifold. Hence, $I_{u,M}'(1)=0$ if and only if
$u\in \mathcal{N}_{\lambda,M}$. We have,\\
$$
\mathcal{N}_{\lambda,M}^{+}=\{u\in
\mathcal{N}_{\lambda,M}:\,I_{u,M}''(1)>0\};
$$
$$
\mathcal{N}_{\lambda,M}^{0}=\{u\in
\mathcal{N}_{\lambda,M}:\,I_{u,M}''(1)=0\};
$$
$$
\mathcal{N}_{\lambda,M}^{-}=\{u\in
\mathcal{N}_{\lambda,M}:\,I_{u,M}''(1)<0\}.
$$
Thus, for each $u\in \mathcal{N}_{\lambda,M}$,
\begin{align*}
I_{u,M}''(1) &=
(p-1)\,M(\|u\|_{{W}^{1,p}}^{p})\,\|u\|_{{W}^{1,p}}^{p}
+p\,M'(\|u\|_{{W}^{1,p}}^{p})\,\|u\|_{{W}^{1,p}}^{2p}
-\lambda\,(q-1)\,\int_{\Omega}f|u|^{q}\,dx\\&-(r-1)\,\int_{\Omega}g|u|^{r}\,dx
\end{align*}
\begin{equation}
\hspace* {1.3cm}
=p\,M'(\|u\|_{{W}^{1,p}}^{p})\,\|u\|_{{W}^{1,p}}^{2p}+\,(p-q)\,M(\|u\|_{{W}^{1,p}}^{p})\,\|u\|_{{W}^{1,p}}^{p}
-\,(r-q)\,\int_{\Omega}g|u|^{r}\,dx
\end{equation}
\begin{equation}
\hspace*
{1.3cm}=p\,M'(\|u\|_{{W}^{1,p}}^{p})\,\|u\|_{{W}^{1,p}}^{2p}
-\,(r-p)\,M(\|u\|_{{W}^{1,p}}^{p})\,\|u\|_{{W}^{1,p}}^{p}+\lambda\,(r-q)\,\int_{\Omega}f|u|^{q}\,dx.
\end{equation}
Define
\begin{equation}
\psi_{\lambda,M}(u)=\langle \,\mathcal{J}_{\lambda,M }'(u)
,u\rangle=M(\|u\|_{{W}^{1,p}}^{p})\,\|u\|_{{W}^{1,p}}^{p}-\lambda\,\int_{\Omega}f|u|^{q}\,dx
-\int_{\Omega}g|u|^{r}\,dx.
\end{equation}
Then for $u\in \mathcal{N}_{\lambda,M}$,
\begin{align*}
\langle\psi_{\lambda,M }'(u)
,u\rangle&=(p-1)\,M(\|u\|_{{W}^{1,p}}^{p})\,\|u\|_{{W}^{1,p}}^{p}+\,p\,M'(\|u\|_{{W}^{1,p}}^{p})\,\|u\|_{{W}^{1,p}}^{2p}
-\lambda\,(q-1)\,\int_{\Omega}f|u|^{q}\,dx\\&-(r-1)\,\int_{\Omega}g|u|^{r}\,dx
\end{align*}
\begin{equation}
\hspace*{-10cm}=I_{u,M}''(1).\\\\
\end{equation}
Also, as proved in Binding, Drabek and Huang \cite {bdh} or in Brown
and
Zhang \cite {kjb-yz}, we have the following lemma.\\\\
{\bf Lemma 2.1.} Suppose that $u_{0}$ is a local minimizer for
$\mathcal{J}_{\lambda,M}$ on $\mathcal{N}_{\lambda,M}$ and that
$u_{0}\notin \mathcal{N}_{\lambda,M}^{0}$. Then
$\mathcal{J}_{\lambda,M }'(u_{0})=0$ in
$W^{-1,p'}(\Omega).$\\\\

Let$\lambda_{0}(a)=\max\{\frac{q\,\lambda_{1}(a)}{p^{\frac{2p-1}{p}}},\frac{q\,\lambda_{2}}{p}\},$
where $\lambda_{1}(a)$ and $\lambda_{2}$are given by\\
$$\lambda_{1}(a)=\frac{p\,S_{q}^{\frac{q}{p}}\,\sqrt[p]{a\,b^{p-1}\,(r-p^{2})\,(r-p)^{p-1}}}{(r-q)\,\|f^{+}\|_{\infty}}
\big(\frac{p\,S_{r}^{\frac{r}{p}}\,\sqrt[p]{a\,b^{p-1}\,(p^{2}-q)\,(p-q)^{p-1}}}{(r-q)\|g^{+}\|_{\infty}}\Big)^{\frac{2p-1-q}{r-2p+1}}$$\\
and\\
$$\lambda_{2}=\frac{b\,S_{q}^{\frac{q}{p}}(r-p)}{(r-q)\,\|f^{+}\|_{\infty}}\,\Big(
\frac{b\,S_{r}^{\frac{r}{p}}\,(p-q)}{(r-q)\,\|g^{+}\|_{\infty}}\Big)^{\frac{p-q}{r-p}},$$\\
then we will state the main theorems.\\\\
{\bf Theorem 2.1.} Suppose that $N$ is any one of $1,p,2p-1$ and
that $r>p^{2}.$ Then for each $a<0$ and $0< \lambda<\lambda_{0}(a)$,
Eq. $(1)$ has two positive solutions $u_{\lambda,M}^{+}\,\in
\mathcal{N}_{\lambda,M}^{+}$ and
$u_{\lambda,M}^{-}\,\in \mathcal{N}_{\lambda,M}^{-}.$\\\\
\hspace{0.6 cm} We define\\
\begin{equation}
\Lambda=\inf\{\|u\|_{W^{1,p}}^{p^{2}}:\,u\in
W_{0}^{1,p}(\Omega),\int_{\Omega}g|u|^{p^{2}}\,dx=1\}.
\end{equation}
then $\Lambda>0$ is achieved by some $\phi_{\Lambda}\in
W_{0}^{1,p}(\Omega)$ with
$\int_{\Omega}g|\phi_{\Lambda}|^{p^{2}}\,dx=1$ and
$\phi_{\Lambda}>0$ a.e. in $\Omega$ according to the compactness of
Sobolev embedding from $W_{0}^{1,p}(\Omega)$ into
$L^{p^{2}}(\Omega)$ and Fatou's
lemma. So,\\
\begin{equation}
\Lambda\int_{\Omega}g|u|^{p^{2}}\,dx\leq\,
\|u\|_{W^{1,p}}^{p^{2}}\,\,\text {for all}\,\,u \in
W_{0}^{1,p}(\Omega),
\end{equation}
and\\
\begin{equation}
\left\{\begin{array}{ll}-\|u\|_{W^{1,p}}^{p}\,\Delta_{p}u=\mu
\|u\|_{W^{1,p}}^{p}\,u,  & x \in \Omega,\\
u=0,  & x\in
\partial \Omega.
\end{array}\right.
\end{equation}
where $\mu$ is an eigenvalue of
Eq. $(10)$, $u\in W_{0}^{1,p}(\Omega)$ is nonzero and eigenvector corresponding to eigenvalue $\mu$ such that\\
$$
\|u\|_{W^{1,p}}^{p}\,\int_{\Omega}|\nabla u|^{p-2}\,\nabla u
\,\nabla
(u.\varphi)\,dx=\mu\,\int_{\Omega}|u|^{p}\,u\,\varphi\,dx,\,\,
\text {for all}\,\, \varphi\in W_{0}^{1,p}(\Omega),
$$
we writhe
$$
I(u)=\|u\|_{W^{1,p}}^{p^{2}},\,\,\text {for}\,\,u\in
\mathbb{S}_{W}=\{u\in
W_{0}^{1,p}(\Omega):\,\int_{\Omega}|u|^{p^{2}}\,dx=1\},
$$
and all distinct eigenvalues of Eq. (10) denoted by
$0<\mu_{1}<\mu_{2}<...,$ we have
$$
\mu_{1}=\inf_{u\in \mathbb{S}_{W^{1,p}}}\,I(u)>0,
$$
where $\mu_{1}$ is simple, isolated and can be achieved at some
$\psi_{1}\in \mathbb{S}_{W}$ and $\psi_{1}>0$ in $\Omega$ ( see
\cite {z-p}).\\\\
Let\\
$$
\hat{\lambda}_{0}(a)=\frac{p\,b\,S_{q}^{\frac{q}{p}}}{(p^{2}-q)\,
\|f^{+}\|_{\infty}}\,\Big(\frac{b\,\Lambda\,(p-q)}{(1-a\lambda)(p^{2}-q)}\Big)^{\frac{p-q}{p}}
$$
Then we have the following result.\\\\
{\bf Theorem 2.2.} Suppose that $N=1,p,2p-1$ and $r=p^{2}.$ Then
\begin{itemize}
\item{(i)} for each $a\geq \frac{1}{\Lambda}$ and $\lambda>0,\mathcal{N}_{\lambda,M}^{+}=\mathcal{N}_{\lambda,M}$ and Eq. $(1)$ has at least one positive solution;
\item{(ii)} for each $a<\frac{1}{\Lambda}$
and $0<\lambda<\frac{1}{p}\,\hat{\lambda}_{0}(a)$, Eq. $(1)$ has two
positive solutions $u_{\lambda,M}^{+}\in
\mathcal{N}_{\lambda,M}^{+}$ and $u_{\lambda,M}^{-}\in
\mathcal{N}_{\lambda,M}^{-},$\,\,\text{ and}\,\, $\lim_{a\to
\frac{1}{\Lambda}^{-}}\,\inf_{u\in
\mathcal{N}_{\lambda,M}^{-}}\,\mathcal{J}_{\lambda,M
}(u)=\infty.$\\
\end{itemize}

For the next result, we note that if $g\geq 0,$ then using Lemma 2
from Alves et al. \cite{ACM}, there exists $C_{*}>0$ independent of
M and $\lambda$ such that
\begin{equation}
\frac{\|u\|_{{W}^{1,p}}^{p}}{(\lambda\,C_{*}^{q}\,\|f^{+}\|_{\infty}+C_{*}^{r}\,\|g^{+}\|_{\infty})\,|\Omega|}\leq
\max\{M(\|u\|_{W^{1,p}}^{p})^{\frac{(p-r+q)}{(r-p)}},M(\|u\|_{W^{1,p}}^{p})^{\frac{p}{(r-p)}}\}.
\end{equation}
Let
$$
L(\lambda)=(\lambda\,C_{*}^{q}\,\|f^{+}\|_{\infty}+C_{*}^{r}\,\|g^{+}\|_{\infty})\,|\Omega|,
$$
$$
\hat{A}_{0}=\frac{b\,(p-q)\,(r-p)}{p(p^{2}-q)}(\frac{p\,S_{r}
^{\frac{r}{p}}}{b\,(p-q)\,(r-q)\,(p^{2}-r)\,\|g^{+}\|_{\infty}})^{\frac{p^{2}-r}{r-p}},
$$
and
$$
A_{0}=\max\{(\frac{b(2r-p)}{r})^{\frac{(p-r+q)}{(r-p)}},(\frac{b\,r}{p})^{\frac{(p-r+q)}{(r-p)}},(\frac{br}{p})^{\frac{p}{(r-p)}}\}.\\\\
$$
Then we have:\\\\
{\bf Theorem 2.3.} Suppose that $r<p^{2}$. Then
\begin{itemize}
\item{(i)} for each $a,\lambda>0$, Eq. $(1)$ has a
positive solution $u_{a,\lambda}$. And, for each $a>\hat{A }_{0}$
and $\lambda>0$,$u_{a,\lambda}\in
\mathcal{N}_{\lambda,M}^{+}=\mathcal{N}_{\lambda,M}$;
\item{(ii)} if $g\geq 0$,then for each $\theta>0$ and
$0<a<\frac{b(r-p)}{rA_{0}L(\theta)}$ there exists
$\tilde{\lambda}_{0}\in(0,\theta]$ such that for
$0<\lambda<\tilde{\lambda}_{0}$, Eq. $(1)$ has two positive
$u_{\lambda,M}^{+}$ and $u_{\lambda,M}^{-}$such that
$u_{\lambda,M}^{\pm}\,\in \mathcal{N}_{\lambda,M}^{\pm}$and
$\|u_{\lambda,M}^{\pm}\|_{W^{1,p}}^{p}<\frac{b(r-p)}{pa}$.\\\\
\end{itemize}
Finally, let
$$
A_{*}=\frac{p^{\frac{r}{(p-r)}}(r-p)^{p}}{r\mathbb{S}}\,\big(\frac{2p-r}{b}\Big)^{\frac{(p^{2}-r)}{(r-p)}},
$$
and\\
$$
\hat{\Lambda}=a\,\big(\frac{b(r-p)}{a(p^{2}-r)}\Big)^{\frac{(p^{2}-q)}{p}}\,\|f\|_{\infty}^{-1}\,S_{q}^{\frac{q}{p}},
$$
where $\mathbb{S}>0$ given by $(36)$, we state our last theorem.\\\\
{\bf Theorem 2.4.} Suppose that $r<p^{2}$ and $f,g>0$. Then for each
$\theta>0$and $0<a<\{\frac{b(r-p)}{rA_{0}L(\theta)},A_{*}\}$ there
exists a positive number $\tilde{\lambda}_{*}\leq\,
\min\{\theta,\hat{\Lambda}\}$such that for
$0<\lambda<\tilde{\lambda}_{*}$, Eq. $(1)$ has three positive
solutions $u_{\lambda,M}^{(1),+}$,$u_{\lambda,M}^{(p),+}$ and
$u_{\lambda,M}^{-}$ such that $u_{\lambda,M}^{(i),+}\in
\mathcal{N}_{\lambda,M}^{+}$,$u_{\lambda,M}^{-}\in
N_{\lambda,M}^{-}$ and
$$
\|u_{\lambda,M}^{(1),+}\|_{W^{1,p}}^{p}<(\frac{p\lambda\,
(r-q)\|f\|_{\infty}S_{q}^{\frac{-q}{p}}}{b(r-p)^{p}})^{\frac{p}{(p-q)}},
$$
$$
S_{r}^{\frac{r}{p(r-p)}}(\frac{pb(r-1)(p-q)}{r(r-q)\|g^{+}\|_{\infty}})^{\frac{1}{(r-p)}}<\|u_{\lambda,M}^{-}\|_{W^{1,p}}^{p}<\frac{b(r-p)}{pa}<\|u_{\lambda,M}^{(p),+}\|_{W^{1,p}}^{p}.
$$
\section {Preliminary Results}
\hspace{0.6 cm}The sequence $\{u_{n}\}$ is a Palais-Smale sequence
for $\mathcal{J}_{\lambda,M}$ on
$W_{0}^{1,p}(\Omega)$ if\\
$$
\mathcal{J}_{\lambda,M}(u_{n})\,\, \text {is bounded and}\,\,
 \mathcal{J}_{\lambda,M}'(u_{n})= o(1) \,\,\text{in} \,\,W^{-1,p'}(\Omega).
$$
Furthermore, if every Palais-Smale sequence for
$\mathcal{J}_{\lambda,M}$ on $W_{0}^{1,p}(\Omega)$ has a strongly
convergent subsequence, then $\mathcal{J}_{\lambda,M}$
satisfies the Palais-Smale condition. Now we have the following results.\\\\
{\bf Lemma 3.1.} Suppose that $M(s)\geq \,m_{0}\,\,\text {for
all}\,\,
 s\geq 0$ and for some
 $m_{0}>0.$ Then
each bounded Palais-Smale sequence for $\mathcal{J}_{\lambda,M}$
on
$W_{0}^{1,p}(\Omega)$ has a strongly convergent subsequence.\\\\
{\bf Proof.} Let $\{u_{n}\}$ be a bounded Palais-Smale sequence
for $\mathcal{J}_{\lambda,M}$ on $W_{0}^{1,p}(\Omega)$. Then by
the compact embedding theorem, there exist subsequences
$\{u_{n}\}$ and $u_{0}\in W_{0}^{1,p}(\Omega)$ such that
$$
u_{n}\rightharpoonup u_{0}\,\,\,\text {Weakly in}\,\,
  W_{0}^{1,p}(\Omega)
$$
and
$$
u_{n}\to u_{0}\,\,\, \text {strongly in} \,\, L^{z}(\Omega)\,\,\text
{for}\,\, 1<z<p^{*}.
$$
then
$$
\int_{\Omega}(\lambda\,
f|u_{n}|^{q-2}\,u_{n}+g|u_{n}|^{r-2}\,u_{n})\,(u_{n}-u_{0})\,dx\to
0,
$$
and since
$$\mathcal{J}_{\lambda,M}'(u_{n}^{p-1})(u_{n}-u_{0})\to 0,$$
We have
$$
M(u_{n})\int_{\Omega} \nabla u_{n}^{p-1}\nabla (u_{n}-u_{0})\to 0.
$$
Thus, $u_{n}\to u_{0}\,\,\,\text {strongly in}\,\,
 W_{0}^{1,p}(\Omega)$.
This completes the proof.\hspace* {.5 cm}$\Box$\\\\
{\bf Lemma 3.2.} Suppose that $M(s)=as^{p-1}+b$. Then,
\begin{itemize}
\item{\bf (i)} If $r\geq p^{2}$, then energy
functional $\mathcal{J}_{\lambda,M}$ is coercive and below on
$\mathcal{N}_{\lambda,M}$;
\item{\bf (ii)} If $r< p^{2}$, then energy
functional $\mathcal{J}_{\lambda,M}$ is coercive and bounded below
on $W_{0}^{1,p}(\Omega)$.
\end{itemize}
{\bf Proof.} $(i)$ For $u\in \mathcal{N}_{\lambda,M}$, we have
$M(\|u\|_{W^{1,p}}^{p})\|u\|_{W^{1,p}}^{p}=\lambda\int_{\Omega}f|u|^{q}\,dx+\int_{\Omega}g|u|^{r}\,dx$.
By the Sobolev inequality,\\
\begin{eqnarray*}
\mathcal{J}_{\lambda,M}(u)&=&\frac{1}{p}\hat{M}(\|u\|_{W^{1,p}}^{p})-\frac{\lambda}{q}\,\int_{\Omega}f|u|^{q}\,dx-\frac{1}{r}\,\int_{\Omega}g|u|^{r}\,dx\\
&=&\frac{1}{p}\hat{M}(\|u\|_{W^{1,p}}^{p})-\frac{1}{p}\,M(\|u\|_{W^{1,p}}^{p})\,\|u\|_{W^{1,p}}^{p}-\lambda\,(\frac{r-q}{rq})\,\int_{\Omega}f|u|^{q}\,dx\\
&\geq&\frac{\|u\|_{W^{1,p}}^{p}}{rp}(\frac{a(r-p^{2})}{p}\,\|u\|_{W^{1,p}}^{p^{2}-p}+b(r-p))-\lambda(\frac{r-q}{rq})\|f^{+}\|_{\infty}
S_{q}^{\frac{-q}{p}}\|u\|_{W^{1,p}}^{q}.
\end{eqnarray*}
and,
\[
\mathcal{J}_{\lambda}(u)\,\geq\,\frac{b(r-p)}{rp}\|u\|_{W^{1,p}}^{p}-\lambda(\frac{r-q}{rq})\|f^{+}\|_{\infty}\,S_{q}^{\frac{-q}{p}}\|u\|_{W^{1,p}}^{q}.
\]
Thus, $\mathcal{J}_{\lambda,M}$ is coercive and bounded below on
$\mathcal{N}_{\lambda,M}$.\\\\
$(ii)$ We have,\\
\begin{eqnarray*}
\mathcal{J}_{\lambda,M}(u)&=&\frac{1}{p}\hat{M}(\|u\|_{W^{1,p}}^{p})-\frac{\lambda}{q}\,\int_{\Omega}f|u|^{q}\,dx-\frac{1}{r}\,\int_{\Omega}g|u|^{r}\,dx\\
&\geq&\,\frac{a}{p^{2}}\|u\|_{W^{1,p}}^{p^{2}}+\frac{b}{p}\,\|u\|_{W^{1,p}}^{p}-\frac{\lambda\|f^{+}\|_{\infty}\,S_{q}^{\frac{-q}{p}}}{q}\|u\|_{W^{1,p}}^{q}-\frac{\|g^{+}\|_{\infty}\,S_{r}^{\frac{-r}{p}}}{q}\,\|u\|_{W^{1,p}}^{r}\\
&=&(\frac{a}{p^{2}}\|u\|_{W^{1,p}}^{p^{2}-r}\,-\frac{S_{r}^{\frac{-r}{p}}\|g^{+}\|_{\infty}}{r})\|u\|_{W^{1,p}}^{r}+(\frac{b}{p}\|u\|_{W^{1,p}}^{p-q}-\frac{\lambda\|f^{+}\|_{\infty}\,S_{q}^{\frac{-q}{p}}}{q})\|u\|_{W^{1,p}}^{q},
\end{eqnarray*}
for all $u\in W_{0}^{1,p}(\Omega)\backslash\{0\}$. Thus,
$\mathcal{J}_{\lambda,M}$ is coercive and bounded below on
$W_{0}^{1,p}(\Omega)$.$\hspace* {.2cm}\Box$\\\\
{\bf Lemma 3.3.} Suppose that $M(s)=as^{p-1}+b$. Then we have
\begin{itemize}
\item{\bf (i)} If $r>p^{2}$ and
$0<\lambda<\max\{\lambda_{1}(a),\lambda_{2}\}$, then for all
$a>0$,\,\,\, $\mathcal{N}_{\lambda,M}^{0}=\emptyset$.
\item{\bf (ii)} if $r=p^{2}$ and $a\geq\frac{1}{\Lambda}$, then  $\text {for all} \lambda>0$,\,\,\,$\mathcal{N}_{\lambda,M}^{+}=\mathcal{N}_{\lambda,M}$.
\item{\bf (iii)} if $r=p^{2}$,$a<\frac{1}{\Lambda}$ and
$0<\lambda<\hat{\lambda}_{0}(a)$, then
$\mathcal{N}_{\lambda,M}^{0}=\emptyset$.
\item{\bf (iv)} if $r<p^{2}$ and $a>\hat{A}_{0}$, then $\text {for
all}\,\,
 \lambda>0$,\,\,\,
$\mathcal{N}_{\lambda,M}^{+}=\mathcal{N}_{\lambda,M}$.\\\\
\end{itemize}

{\bf Proof.} $(i)$ If $r>p^{2}$ and $u\in
\mathcal{N}_{\lambda,M}^{0}$, then
\begin{equation}
\begin{aligned}
(r-q)\|g^{+}\|_{\infty}S_{r}^{\frac{-r}{p}}\|u\|_{W^{1,p}}^{r} &\geq
a(p^{2}-q)\|u\|_{W^{1,p}}^{p^{2}}+b(p-q)(p-1)\|u\|_{W^{1,p}}^{p}\\
&\geq
\begin{cases}
p\sqrt[p]{ab^{p-1}(p^{2}-q)(p-q)^{p-1}}\|u\|_{W^{1,p}}^{2p-1},\\b(r-p)\|u\|_{W^{1,p}}^{p},\\
\end{cases}
\end{aligned}
\end{equation}
and
\begin{equation}
\begin{aligned}
\lambda(r-q)\|f^{+}\|_{\infty}S_{q}^{\frac{-q}{p}}\|u\|_{W^{1,p}}^{q}&\geq
a(r-p^{2})\|u\|_{W^{1,p}}^{p^{2}}+b(r-p)(p-1)\|u\|_{W^{1,p}}^{p}\\
&\geq
\begin{cases}
p\sqrt[p]{ab^{p-1}(r-p^{2})(r-p)^{p-1}}\|u\|_{W^{1,p}}^{2p-1},\\b(r-p)\|u\|_{W^{1,p}}^{p}.\\
\end{cases}
\end{aligned}
\end{equation}
By $(12)$ and $(13)$\,\,\,for all $u\in
\mathcal{N}_{\lambda,M}^{0}$, we
have\\
\begin{eqnarray*}
(\frac{p
S_{r}^{\frac{r}{p}}\,\sqrt[p]{ab^{p-1}(p^{2}-q)(p-q)^{p-1}}}{(r-q)\|g^{+}\|_{\infty}})^{\frac{1}{r-2p+1}}\leq\|u\|_{W^{1,p}}\leq\,(\frac{\lambda(r-q)\,\|f^{+}\|_{\infty}}{p
S_{q}^{\frac{q}{p}}\,\sqrt[p]{ab^{p-1}(r-p^{2})\,(r-p)^{p-1}}})^{\frac{1}{2p-q-1}}
\end{eqnarray*}
and
\begin{eqnarray*}
(\frac{bS_{r}^{\frac{r}{p}}(p-q)}{(r-q)\|g^{+}\|_{\infty}})
^{\frac{1}{(r-p)}}\leq\|u\|_{W^{1,p}}\leq(\frac{\lambda(r-q)\|f^{+}\|_{\infty}}{S_{q}^{\frac{-q}{p}}b(r-p)})^{\frac{1}{p-q}}.
\end{eqnarray*}
Hence, if $\mathcal{N}_{\lambda,M}^{0}$ is nonempty, then the
inequality
$\lambda\geq\,\max\{\lambda_{1}(a),\lambda_{2}\}$ must hold.\\
$(ii)$ If $r=p^{2}$ and $a\,\geq\,\frac{1}{\Lambda}$, then for all
$u\in \mathcal{N}_{\lambda,M}$ we have\\
\begin{eqnarray*}
I_{\lambda,M}''(1)&=&a\,(p^{2}-q)\,\|u\|_{W^{1,p}}^{p^{2}}+b(p-q)\,\|u\|_{W^{1,p}}^{p}-(p^{2}-q)\,\int_{\Omega}g|u|^{r}\,dx\\
&\geq&\frac{(a\Lambda-1)(p^{2}-q)}{\Lambda}\,\|u\|_{W^{1,p}}^{p^{2}}+b(p-q)\,\|u\|_{W^{1,p}}^{p}>0,
\end{eqnarray*}
Thus, $\mathcal{N}_{\lambda,M}^{+}=\mathcal{N}_{\lambda,M}$ for all $\lambda>0$.\\
$(iii)$ If $r=p^{2}$, $a<\frac{1}{\Lambda}$ and $u\in
\mathcal{N}_{\lambda,M}^{0}$, then
\begin{eqnarray*}
b(p-q)\,\|u\|_{W^{1,p}}^{p}=(p^{2}-q)(\int_{\Omega}g|u|^{r}\,dx-a\,\|u\|_{W^{1,p}}^{p^{2}})
\end{eqnarray*}
\begin{equation}
\hspace* {1 cm}
\leq\,\frac{(1-a\Lambda)(p^{2}-q)}{\Lambda}\,\|u\|_{W^{1,p}}^{p^{2}},
\end{equation}
and
\begin{equation}
\|u\|_{W^{1,p}}^{p}\leq\frac{\lambda(p^{2}-q)\|f^{+}\|_{\infty}}{b(p^{2}-p)S_{q}^{\frac{q}{p}}}\,\|u\|_{W^{1,p}}^{q}.
\end{equation}
By $(14)$ and $(15)$ for all $u\in \mathcal{N}_{\lambda,M}^{0}$, we
have
\begin{eqnarray*}
(\frac{b\,\Lambda\,(p-q)}{(1-a\Lambda)(p^{2}-q)})^{\frac{1}{p^{2}-p}}\,\leq\|u\|_{W^{1,p}}\,\leq\,(\frac{\lambda(p^{2}-q)\|f^{+}\|_{\infty}}{b(p^{2}-p)S_{q}^{\frac{q}{p}}})^{\frac{1}{p-q}}
\end{eqnarray*}
Hence, if $\mathcal{N}_{\lambda,M}^{0}$ is nonempty, then the
inequality
$\lambda\geq\hat{\lambda}_{0}(a)$ must be hold.\\
$(iv)$ If $r<p^{2}$ and $u\in \mathcal{N}_{\lambda,M}$, then
\begin{eqnarray*}
I_{\lambda,M}''(1)&=&a(p^{2}-q)\,\|u\|_{W^{1,p}}^{p^{2}}+b(p-q)\,\|u\|_{W^{1,p}}^{p}-(r-q)\,\int_{\Omega}g|u|^{r}\,dx\\
&\geq&a(p^{2}-q)\,\|u\|_{W^{1,p}}^{p^{2}}+b(p-q)\,\|u\|_{W^{1,p}}^{p}-(r-q)\|g^{+}\|_{\infty}\,S_{r}^{\frac{-r}{p}}\|u\|_{W^{1,p}}^{r}\\
&=&\|u\|_{W^{1,p}}^{p}[a(p^{2}-q)\,\|u\|_{W^{1,p}}^{p^{2}-p}+b(p-q)-(r-q)\|g^{+}\|_{\infty}\,S_{r}^{\frac{-r}{p}}\|u\|_{W^{1,p}}^{r-p}].
\end{eqnarray*}
Hence, if $a>\hat{A}_{0}$, then we have
\begin{eqnarray*}
a(p^{2}-q)\,\|u\|_{W^{1,p}}^{p^{2}-p}+b(p-q)-(r-q)\,\|g^{+}\|_{\infty}\,S_{r}^{\frac{-r}{p}}\|u\|_{W^{1,p}}^{r-p}\,\,\,\,\text{for
all}\,\,  u\in \mathcal{N}_{\lambda,M},
\end{eqnarray*}
and so
$$
I_{\lambda,M}''(1)>0 \,\,\text {for all} \,\,a>A_{0}
\,\,\text{and}\,\, u\in \mathcal{N}_{\lambda,M}.
$$
Thus, $\mathcal{N}_{\lambda,M}^{+}=\mathcal{N}_{\lambda,M}$for
all $\lambda>0$.
\section{Non-emptiness of submanifolds}
\hspace{0.6 cm}Now we state the results for
$\mathcal{N}_{\lambda,M}^{0}$,$\mathcal{N}_{\lambda,M}^{-}$ and
$\mathcal{N}_{\lambda,M}^{+}$ that are non-empty for various $r,a$
and $\lambda$.\\\\
{\bf Theorem 4.1.}
\begin{itemize}
\item{\bf (i)} If $r>p^{2}$ and
$0<\lambda<\max\{\lambda_{1}(a),\lambda_{2}\}$, then
$\mathcal{N}_{\lambda,M}=\mathcal{N}_{\lambda,M}^{+}\cup
\mathcal{N}_{\lambda,M}^{-}$ and
$\mathcal{N}_{\lambda,M}^{\pm}\neq\emptyset$ for all $a>0$.
\item{\bf (ii)} If $r=p^{2}$, $a<\frac{1}{\Lambda}$ and
$0<\lambda<\hat{\lambda}_{0}(a)$, then
$\mathcal{N}_{\lambda,M}=\mathcal{N}_{\lambda,M}^{+}\cup
\mathcal{N}_{\lambda,M}^{-}$ and
$\mathcal{N}_{\lambda,M}^{\pm}\neq\emptyset$.
\item{\bf (iii)}If $r=p^{2}$, $a\geq\frac{1}{\Lambda}$, then for all
$\lambda>0$,
$\mathcal{N}_{\lambda,M}^{+}=\mathcal{N}_{\lambda,M}\neq\emptyset$.
\item{\bf (iv)}If $r<p^{2}$ and $a>\hat{A}_{0}$,
then for all $\lambda>0$,
$\mathcal{N}_{\lambda,M}^{+}=\mathcal{N}_{\lambda,M}\neq\emptyset$.\\
\end{itemize}
\hspace*{0.6 cm}To prove Theorem 4.1. (i) we state the following lemmas:\\\\
{\bf Lemma 4.2.} Suppose that $r>p^{2}$ and
$0<\lambda<\max\{\lambda_{1}(a),\lambda_{2}\}$. Then for each $u\in
{W_{0}}^{1,p}(\Omega)$ with $\int_{\Omega}g|u|^{r}\,dx>0$, there
exists $t_{a,\rm max}>0$ such that\\
\begin{itemize}
\item{\bf (i)} If $\int_{\Omega}f|u|^{q}\,dx\,\leq 0$, then there is a
unique $t^{-}>t_{a,{\rm max}}$ such that $t^{-}u\in
\mathcal{N}_{\lambda,M}^{-}$ and
$$
\mathcal{J}_{\lambda,M}(t^{-}{u})=\sup_{t\geq
0}\,\mathcal{J}_{\lambda,M}(tu);
$$
\item{\bf (ii)} If$\int_{\Omega}f|u|^{q}\,dx>0$, then there are unique
$t^{+}$ and $t^{-}$ with $0<t^{+}<t_{a,\max}<t^{-}$such that
$t^{\pm}{u}\in \mathcal{N}_{\lambda,M}^{\pm}$ \,and\,\\
$$
\mathcal{J}_{\lambda,M}(t^{+}{u})=\inf_{0\leq t\leq
t_{a,\max}}\,\mathcal{J}_{\lambda,M}(tu);\,\,\text {and}\,\,
\mathcal{J}_{\lambda,M}(t^{-}u)=\sup_{t\geq
t_{a,max}}\mathcal{J}_{\lambda,M}(tu).
$$
\end{itemize}
{\bf Proof.}
Case (A): $\lambda_{0}=\lambda_{2}$.\\
\hspace{0.6 cm} Fix $u\in {W_{0}}^{1,p}(\Omega)$ with
$\int_{\Omega}g|u|^{r}\,dx>0$. Let
\begin{eqnarray*}
h_{a}(t)=at^{p^{2}-q}\|u\|_{W^{1,p}}^{p^{2}}+bt^{p-q}\,\|u\|_{W^{1,p}}^{p}-t^{r-q}\,\int_{\Omega}g|u|^{r}\,dx\,\,\text
  {for}   a,t\geq0.
\end{eqnarray*}
Clearly, $tu\in \mathcal{N}_{\lambda,M}$ if and only if
$h_{a}(t)=\lambda\,\int_{\Omega}f|u|^{q}\,dx$. We have $h_{a}(0)=0$
and $h_{a}(t)\to -\infty$ as $t\to \infty$. Since
$\int_{\Omega}g|u|^{r}\,dx>0$, $r>p^{2}$ and
\begin{eqnarray*}
h_{a}'(t)=t^{p-q-1}\Big(a\,(p^{2}-q)t^{p^{2}-p}\,\|u\|_{W^{1,p}}^{p^{2}}+b(p-q)\,\|u\|_{W^{1,p}}^{p}-(r-q)\,t^{r-p}\,\int_{\Omega}g|u|^{r}\,dx\Big),
\end{eqnarray*}
there is a unique $t_{a,\max}>0$ such that $h_{a}(t)$ achieves its
maximum at $t_{a,\max}$, increasing for $t\in [0,t_{a,\max})$ and
decreasing for $t\in (t_{a,\max},\infty)$ with $\lim_{t\to \infty}\,
h_{a}(t)=-\infty$. Clearly, if $tu\in \mathcal{N}_{\lambda,M}$, then
$t^{q-1}\,h_{a}'(t)=I_{u,M}''(t)$. Hence $tu\in
\mathcal{N}_{\lambda,M}^{+}$(or$\mathcal{N}_{\lambda,M}^{-}$)if and
only if $h_{a}'(t)>0$ (or $<0$). Moreover,
$$
t_{0,\max}=(\frac{b(p-q)\|u\|_{W^{1,p}}^{p}}{(r-q)\int_{\Omega}g|u|^{r}\,dx})^{\frac{1}{(r-p)}},
$$
and
\begin{align*}
h_{0}(t_{0,\max})&= b\Big(\frac{b(p-q)\|u\|_{W^{1,p}}^{p}}{(r-q)\int_{\Omega}g|u|^{r}\,dx}\Big)^{\frac{(p-q)}{(r-p)}}\,\|u\|_{W^{1,p}}^{p}-(\frac{b(p-q)\|u\|_{W^{1,p}}^{p}}{(r-q)\,\int_{\Omega}g|u|^{r}\,dx})^{\frac{r-q}{(r-p)}}\,\int_{\Omega}g|u|^{r}\,dx\\
&=\|u\|_{W^{1,p}}^{q}\,[(\frac{p-q}{r-q})^{\frac{p-q}{r-p}}-(\frac{p-q}{r-q})^{\frac{r-q}{r-p}}]\,(\frac{\|u\|_{W^{1,p}}^{r}}{\int_{\Omega}g|u|^{r}\,dx})^{\frac{p-q}{r-p}}\,b^{\frac{(r-q)}{(r-p)}}
\end{align*}
\begin{equation}
\hspace* {-5
cm}\geq\|u\|_{W^{1,p}}^{q}\frac{b(r-p)}{(r-q)}(\frac{bS_{r}^{\frac{r}{p}}(p-q)}{(r-q)\|g^{+}\|_{\infty}})^{\frac{(p-q)}{(r-p)}}.
\end{equation}
Case (A-i): $\int_{\Omega}f|u|^{q}\,dx\,\leq 0$.\\
\hspace{0.6 cm} There is a unique $t^{-}>t_{a,max}$ such that
$h_{a}(t^{-})=\lambda\int_{\Omega}f|u|^{q}\,dx$ and
$h_{a}'(t^{-})<0$. Now,
\begin{eqnarray*}
I_{t^{-}u,M}'(1)&=&t^{-}I_{u,M}'(t^{-})\\&=&M(\|t^{-}u\|_{W^{1,p}}^{p})\|t^{-}u\|_{W^{1,p}}^{p}-\lambda\,\int_{\Omega}f|t^{-}u|^{q}\,dx-\int_{\Omega}g|t^{-}u|^{r}\,dx\\
&=&(t^{-})^{q}[h_{a}(t^{-})-\lambda\int_{\Omega}f|t^{-}u|^{q}\,dx]=0,
\end{eqnarray*}
and
\begin{eqnarray*}
I_{t^{-}u,M}''(1)=(t^{-})^{2}I_{u,M}'(t^{-})=(t^{-})^{q+1}h_{a}'(t^{-})<0.
\end{eqnarray*}
Thus $t^{-}u\in \mathcal{N}_{\lambda,M}^{-}$. Since for
$t>t_{a,\max}$, we have $h_{a}'(t)<0$ and $h_{a}''(t)<0$.
Subsequently,
$$
\mathcal{J}_{\lambda,M}(t^{-}u)=I_{u,M}(t^{-})=\sup_{t\geq
0}I_{u,M}(t)=\sup_{t\geq0}\mathcal{J}_{\lambda,M}(tu).
$$
Case (A-ii): $\int_{\Omega}f|u|^{q}\,dx>0$. By $(16)$ and
\begin{eqnarray*}
h_{a}(0)&=&0<\lambda\int_{\Omega}f|u|^{q}\,dx
\\&\leq&\lambda\|f^{+}\|_{\infty}S_{q}^{\frac{-q}{p}}\|u\|_{W^{1,p}}^{q}\\&<&\|u\|_{W^{1,p}}^{q}\frac{b(r-p)}{(r-q)}(\frac{bS_{r}^{\frac{r}{p}}(p-q)}{(r-q)\|g^{+}\|_{\infty}})^{\frac{(r-q)}{(r-p)}}\\
&\leq&\,h_{0}(t_{0,\max})<h_{a}(t_{a,\max}).
\end{eqnarray*}
there are unique $t^{+}$ and $t^{-}$ such that
$0<t^{+}<t_{a,max}<t^{-}$,
$$
h_{a}(t^{+})=\lambda\int_{\Omega}f|u|^{q}\,dx=h_{a}(t^{-}),
$$
and
$$
h_{a}'(t^{+})>0>h_{a}'(t^{-}).
$$
Similar to the argument in part $(A-i)$, we conclude that
$t^{+}u\in \mathcal{N}_{\lambda,M}^{+}$ and $t^{-}u\in
\mathcal{N}_{\lambda,M}^{-}$. Moreover,
\begin{eqnarray*}
\mathcal{J}_{\lambda,M}(t^{-}u)\geq\mathcal{J}_{\lambda,M}(tu)\geq\mathcal{J}_{\lambda,M}(t^{+}u)\,\,\text
 {for each}\,\,   t\in [t^{+},t^{-}],
\end{eqnarray*}
and
$\mathcal{J}_{\lambda,M}(t^{+}u)\leq\mathcal{J}_{\lambda,M}(tu)$
for each $t\in [0,t^{+}]$. Thus,
\begin{eqnarray*}
\mathcal{J}_{\lambda,M}(t^{+}u)=\inf_{0\leq t\leq
t_{a,\max}}\mathcal{J}_{\lambda,M}(t^{+}u)\,\,\text  {and}\,\,
 \mathcal{J}_{\lambda,M}(t^{-}u)=\sup_{t\geq
t_{a,\max}}\mathcal{J}_{\lambda,M}(tu).
\end{eqnarray*}
Case (B): $\lambda_{0}=\lambda_{1}(a)$.\\
\hspace{0.6 cm}Fix $u\in W_{0}^{1,p}(\Omega)$ with
$\int_{\Omega}g|u|^{p}\,dx>0$. Let
\begin{eqnarray*}
m_{a}(t)=p\sqrt[p]{ab}t^{2p-q-1}\|u\|_{W^{1,p}}^{2p-1}-t^{r-q}\int_{\Omega}g|u|^{r}\,dx\,\,\text
{for} t\geq0.
\end{eqnarray*}
Then, $m_{a}(t)\leq h_{a}(t)$ for all $a>0$ and $t\geq 0$. We have
$m_{a}(0)=0$ and $m_{a}(t)\to -\infty$ as $t\to \infty$. Since
$\int_{\Omega}g|u|^{p}\,dx>0$,$r>p^{2}$ and
$$
m_{a}'(t)=t^{1-q}(p\sqrt[p]{ab}(2p-q-1)t^{2p-3}\|u\|_{W^{1,p}}^{2p-1}-(r-q)t^{r-2}\int_{\Omega}g|u|^{r}\,dx),
$$
there is a unique
$$
\tilde{t}_{a,\max}=(\frac{p\sqrt[p]{ab}(2p-q-1)t^{2p-3}\|u\|_{W^{1,p}}^{2p-1}}{(r-q)\int_{\Omega}g|u|^{r}\,dx})^{\frac{1}{r-2p+1}}>0
$$
such that $m_{a}(t)$ achives its maximum at $\tilde{t}_{a,\max}$,
increasing for $t\in [0,\tilde{t}_{a,\max})$ and decreasing for
$t\in (\tilde{t}_{a,\max},\infty)$. Moreover,
\begin{equation}
\begin{aligned}
m_{a}(\tilde{t}_{a,\max})&=p\sqrt[p]{ab}\|u\|_{W^{1,p}}^{q}(\frac{p\sqrt[p]{ab}(2p-q-1)
\|u\|_{W^{1,p}}^{r}}{(r-q)\int_{\Omega}g|u|^{r}\,dx})^{\frac{(2p-q-1)}{(r-2p+1)}}\\
&-p\sqrt[p]{ab}\|u\|_{W^{1,p}}^{q}\,\frac{2p-q-1}{r-q}\,(\frac{p\sqrt[p]{ab}(2p-q-1)\|u\|_{W^{1,p}}^{r}}{(r-q)\int_{\Omega}g|u|^{r}\,dx})^{\frac{(2p-q-1)}{(r-2p+1)}}\\
&=\|u\|_{W^{1,p}}^{q}\,\frac{p(r-2p+1)\sqrt[p]{ab}}{(r-q)}\,(\frac{p\sqrt[p]{ab}(2p-q-1)\|u\|_{W^{1,p}}^{r}}{(r-q)\int_{\Omega}g|u|^{r}\,dx})^{\frac{(2p-q-1)}{(r-2p+1)}}
\\&\geq\|u\|_{W^{1,p}}^{q}\,\frac{p(r-2p+1)\sqrt[p]{ab}}{(r-q)}\,(\frac{p\sqrt[p]{ab}(2p-q-1)S_{r}^{\frac{r}{p}}}{(r-q)\|g^{+}\|_{\infty}})^{\frac{(2p-q-1)}{(r-2p+1)}}.
\end{aligned}
\end{equation}
Case (B-i):$\int_{\Omega}f|u|^{q}\,dx\leq 0$.\\
\hspace{0.6 cm}By $h_{a}(0)=0$ and $h_{a}(t)\to -\infty$ as $t\to
\infty $, there is a unique $t^{-}>t_{a,\max}$ such that
$h_{a}(t^{-})=\lambda\int_{\Omega}f|u|^{q}\,dx$ and
$h_{a}'(t^{-})<0$. Repeating the argument in part $(A-i)$, we have
$t^{-}u\in \mathcal{N}_{\lambda,M}^{-}$ and
$$
\mathcal{J}_{\lambda,M}(t^{-}u)=I_{u,M}(t^{-})=\sup_{t\geq0}I_{u,M}(t)=\sup_{t\geq0}\mathcal{J}_{\lambda,M}(tu).
$$
Case (B-ii):$\int_{\Omega}f|u|^{q}\,dx>0$. By $(16)$ and
\begin{eqnarray*}
h_{a}(0)=0&<&\lambda\int_{\Omega}f|u|^{q}\,dx\leq\lambda\|f^{+}\|_{\infty}S_{q}^{\frac{-q}{p}}\|u\|_{W^{1,p}}^{q}\\
&<&\|u\|_{W^{1,p}}^{q}\frac{p\sqrt[p]{ab^{p-1}(r-p^{2})(r-p)^{p-1}}}{(r-q)}\,(\frac{pS_{r}^{\frac{r}{p}}\sqrt[p]{ab^{p-1}(p^{2}-q)(p-q)^{p-1}}}{(r-q)\|g^{+}\|_{\infty}})^{\frac{(2p-q-1)}{(r-2p+1)}}\\
&\leq&\,\|u\|_{W^{1,p}}^{q}\frac{p\sqrt[p]{ab}(r-2p+1)}{(r-q)}\,(\frac{p(2p-q-1)S_{r}^{\frac{r}{p}}\sqrt[p]{ab}}{(r-q)\|g^{+}\|_{\infty}})^{\frac{(2p-q-1)}{(r-2p+1)}}\\
&\leq&\, m_{a}(\tilde{t}_{a,\max})<h_{a}(t_{a,max}),
\end{eqnarray*}
there are unique $t^{+}$ and $t^{-}$ such that
$0<t^{+}<t_{a,\max}<t^{-}$,
$$
h_{a}(t^{+})=\lambda\int_{\Omega}f|u|^{q}\,dx=h_{a}(t^{-}),
$$
and
$$
h_{a}'(t^{+})>0>h_{a}'(t^{-}).
$$
Repeating the same argument of part $(A-i)$, we conclude that
$t^{+}u\in \mathcal{N}_{\lambda,M}^{+}$ and $t^{-}u\in
\mathcal{N}_{\lambda,M}^{-}$. Moreover,
\begin{eqnarray*}
\mathcal{J}_{\lambda,M}(t^{-}u)\geq\mathcal{J}_{\lambda,M}(tu)\geq\mathcal{J}_{\lambda,M}(t^{+}u)\,\,\text{for
each}\,\,  t\in[t^{+},t^{-}],
\end{eqnarray*}
and $\mathcal{J}_{\lambda,M}(t^{+}u)\leq\mathcal{J}_{\lambda,M}(tu)$
for each $t\in[0,t^{+}]$. Thus,
\begin{eqnarray*}
\mathcal{J}_{\lambda,M}(t^{+}u)=\inf_{0\leq t\,\leq\,t_{a,
\max}}\mathcal{J}_{\lambda,M}(tu)\,\,\text {and}\,\,
\mathcal{J}_{\lambda,M}(t^{-}u)=\sup_{t\geq
t_{a,\max}}\mathcal{J}_{\lambda,M}(tu).
\end{eqnarray*}
This completes the proof.$\hspace* {.2cm}\Box$\\\\
{\bf Lemma 4.3.} Suppose that $r>p^{2}$ and
$0<\lambda<\max\{\lambda_{1}(a),\lambda_{2}\}$. Then for each
$W_{0}^{1,p}(\Omega)$ with $u\in \int_{\Omega}f|u|^{q}\,dx>0$, there
exists $\bar{t}_{a,\max}>0$ such that
\begin{itemize}
\item{\bf (i)} If $\int_{\Omega}g|u|^{r}\,dx\leq 0$ then there is a
unique $0<t^{+}<\bar{t}_{a,\max}$ such that $t^{+}u\in
\mathcal{N}_{\lambda,M}^{+}$ and
$$
\mathcal{J}_{\lambda,M}(t^{+}u)=\inf_{t\geq0}\mathcal{J}_{\lambda,M}(tu);
$$
\item{\bf (ii)} If $\int_{\Omega}g|u|^{r}\,dx>0$ then there are
unique $t^{+}$  and  $t^{-}$ with $0<t^{+}<\bar{t}_{\max}<t^{-}$
such that $t^{\pm} u\in \mathcal{N}_{\lambda,M}^{\pm}$ and
\begin{eqnarray*}
\mathcal{J}_{\lambda,M}(t^{+}u)=\inf_{0\leq
t\leq\,t_{a,max}}\mathcal{J}_{\lambda,M}(tu) \,\,\text {and}\,\,
 \mathcal{J}_{\lambda,M}(t^{-}u)=\sup_{t\geq
t_{a,\max}}\mathcal{J}_{\lambda,M}(tu).
\end{eqnarray*}
\end{itemize}
{\bf Proof.} Fix $u\in W_{0}^{1,p}(\Omega)$ with
$\int_{\Omega}f|u|^{q}\,dx>0$. Let
\begin{eqnarray*}
\bar{h_{a}}(t)=at^{p^{2}-r}\|u\|_{W^{1,p}}^{p^{2}}+bt^{p-r}\|u\|_{W^{1,p}}^{p}-t^{q-r}\lambda\int_{\Omega}f|u|^{q}\,dx\,\,\text
 {for}\,\, t>0\,\,\text {and}  a\geq 0.
\end{eqnarray*}
Clearly, $\bar{h_{a}}(t)\to -\infty$ as $t\to 0^{+}$ and
$\bar{h_{a}}(t)\to 0$ as $t\to \infty$. Since
\begin{eqnarray*}
\bar{h_{a}}'(t)=t^{p-r-1}(a(p^{2}-r)t^{p^{2}-p}\|u\|_{W^{1,p}}^{p^{2}}+b(p-r)\|u\|_{W^{1,p}}^{p}-(q-r)t^{q-p}\lambda\int_{\Omega}f|u|^{q}\,dx,
\end{eqnarray*}
there is a unique $\bar{t}_{a,\max}>0$ such that $\bar{h}_{a}(t)$
achives its maximum at $\bar{t}_{a,\max}$, increasing for $t\in
[0,\bar{t}_{a,\max})$ and decreasing for $t\in
(\bar{t}_{a,\max},\infty)$. Moreover,
$$
\bar{t}_{0,\max}=(\frac{(r-q)\lambda\int_{\Omega}f|u|^{q}\,dx}{b(r-p)\|u\|_{W^{1,p}}^{p}})^{\frac{1}{(p-q)}},
$$
and
\begin{eqnarray*}
\bar{h}_{0}(\bar{t}_{0,\max})&=&b(\frac{b(r-p)\|u\|_{W^{1,p}}^{p}}{(r-q)\lambda\int_{\Omega}f|u|^{q}\,dx})^{\frac{r-p}{p-q}}\,\|u\|_{W^{1,p}}^{p}-(\frac{b(r-p)\|u\|_{W^{1,p}}^{p}}{(r-q)\lambda\int_{\Omega}f|u|^{q}\,dx})^{\frac{r-q}{p-q}}\,\lambda\int_{\Omega}f|u|^{q}\,dx\\
&=&\|u\|_{W^{1,p}}^{r}\,\frac{b(p-q)}{(r-q)}\,(\frac{b(r-p)\|u\|_{W^{1,p}}^{q}}{(r-q)\lambda\int_{\Omega}f|u|^{q}\,dx})^{\frac{r-p}{p-q}}\\
&\geq&\|u\|_{W^{1,p}}^{r}\,\frac{b(p-q)}{(r-q)}\,(\frac{b(r-p)S_{q}^{\frac{q}{p}}}{\lambda\,(r-q)\|f^{+}\|_{\infty}})^{\frac{r-p}{p-q}}.
\end{eqnarray*}
The results of Lemma $4.3$ are obtained by repeating the same
argument of Lemma $4.2$.\\
\hspace{0.6 cm} For the proof of
Theorem 4.1(ii), we require the following two lemmas:\\\\
{\bf Lemma 4.4.} Suppose that $r=p^{2}$, $a<\frac{1}{\lambda}$ and
$0<\lambda<\hat{\lambda}_{0}(a)$. Let $\phi_{\Lambda}>0$ as in
$(8)$. Then, there exists $\hat{t}_{\max}>0$ such that
\begin{itemize}
\item{\bf (i)} If $\int_{\Omega}f|\phi_{\Lambda}|^{q}\,dx\leq 0$ then there is a unique $t^{-}>\hat{t}_{a,\max}$ such that $t^{-}\phi_{\Lambda}\in
\mathcal{N}_{\lambda,M}^{-}$ and
$$
\mathcal{J}_{\lambda,M}(t^{-}\phi_{\Lambda})=\sup_{t\geq0}\mathcal{J}_{\lambda,M}(t\phi_{\Lambda});
$$
\item{\bf (ii)} If $\int_{\Omega}f|\phi_{\Lambda}|^{q}\,dx>0$, then
there are unique $t^{+}$ and $t^{-}$ with
$0<t^{+}<\hat{t}_{\max}<t^{-}$ such that $t^{\pm}\phi_{\Lambda}\in
\mathcal{N}_{\lambda,M}^{\pm}$ and\\
\begin{eqnarray*}
\mathcal{J}_{\lambda,M}(t^{+}\phi_{\Lambda})=\inf_{0\leq t\leq
t_{a,max}}\mathcal{J}_{\lambda,M}(t\phi_{\Lambda})\,\,\text
{and}\,\, \mathcal{J}_{\lambda,M}(t^{-}\phi_{\Lambda})=\sup_{t\geq
t_{a,\max}}\mathcal{J}_{\lambda,M}(t\phi_{\Lambda}).
\end{eqnarray*}
\end{itemize}
{\bf Proof.} Let\\
\begin{eqnarray*}
\hat{h}(t)&=&at^{p^{2}-q}\|\phi_{\Lambda}\|_{W^{1,p}}^{p^{2}}+bt^{p-q}\|\phi_{\Lambda}\|_{W^{1,p}}^{p}-t^{p^{2}-q}\int_{\Omega}g|\phi_{\Lambda}|^{p^{2}}\,dx\\&=&bt^{p-q}\|\phi_{\Lambda}\|_{W^{1,p}}^{p}-t^{p^{2}-q}(\int_{\Omega}g|\phi_{\Lambda}|^{p^{2}}\,dx-a\|\phi_{\Lambda}\|_{W^{1,p}}^{p^{2}})\,\text
{for}  t\geq 0.
\end{eqnarray*}
Then by $(8)$ and $(9)$,
$$
\int_{\Omega}g|\phi_{\Lambda}|^{p^{2}}\,dx-a\|\phi_{\Lambda}\|_{W^{1,p}}^{p^{2}}=1-a\Lambda>0,
$$
we have $\hat{h}(0)=0$ and $\hat{h}(t)\to -\infty$ as $t\to \infty$.
Since\\
\begin{eqnarray*}
\hat{h}'(t)&=&b(p-q)t^{p-q-1}\|\phi_{\Lambda}\|_{W^{1,p}}^{p}-(p^{2}-q)t^{p^{2}-q-1}(\int_{\Omega}g|\phi_{\Lambda}|^{p^{2}}\,dx-a\|\phi_{\Lambda}\|_{W^{1,p}}^{p^{2}})\\
&=&b(p-q)t^{p-q-1}\Lambda^{\frac{1}{p}}-(p^{2}-q)t^{p^{2}-q-1}(1-a\Lambda),
\end{eqnarray*}
there is a unique $\hat{t}_{max}>0$ such that $\hat{h}(t)$ achieves
its maximum at $\hat{t}_{\max}$, increasing for $t\in
[0,\hat{t}_{\max})$ and decreasing for $t\in
(\hat{t}_{{\max},\infty})$. Moreover,
\begin{eqnarray*}
\hat{h}(\hat{t}_{\max})&=&b(\frac{b(p-q)\Lambda^{\frac{1}{p}}}{(p^{2}-q)(1-a\Lambda)})^{\frac{p-q}{p^{2}-q}}\,\Lambda^{\frac{1}{p}}-(1-a\Lambda)(\frac{b(p-q)\Lambda^{\frac{1}{p}}}{(p^{2}-q)(1-a\Lambda)})^{\frac{p^{2}-q}{p^{2}-p}}\\
&=&\Lambda^{\frac{pq-q}{p^{3}-p^{2}}}[\frac{b^{\frac{p^{2}-q}{p^{2}-p}}\,(p-q)^{\frac{p-q}{p^{2}-p}}}{(p^{2}-q)^{\frac{p-q}{p^{2}-p}}}\,(\frac{\Lambda}{1-a\Lambda})^{\frac{p-q}{p^{2}-p}}]\\
&-&\Lambda^{\frac{pq-q}{p^{3}-p^{2}}}\big[\frac{b^{\frac{p^{2}-q}{p^{2}-p}}\,(p-q)^{\frac{p^{2}-q}{p^{2}-p}}}{(p^{2}-q)^{\frac{p^{2}-q}{p^{2}-p}}}\,(\frac{\Lambda}{1-a\Lambda})^{\frac{p-q}{p^{2}-p}}\Big]\\
&=&\Lambda^{\frac{pq-q}{p^{3}-p^{2}}}\,b^{\frac{p^{2}-q}{p^{2}-p}}\,[(\frac{p-q}{p^{2}-q})^{\frac{p-q}{p^{2}-p}}-(\frac{p-q}{2p-q})^{\frac{2p-q}{p}})]\,(\frac{\Lambda}{1-a\Lambda})^{\frac{p-q}{p^{2}-p}}\\
&\geq\,&
\frac{(p^{2}-p)\Lambda^{\frac{q}{p^{2}}}\,b^{\frac{p^{2}-q}{p^{2}-p}}}{p^{2}-q}\,\big(\frac{\Lambda\,(p-q)}{(p^{2}-q)\,(1-a\Lambda)}\Big)^{\frac{p-q}{p^{2}-p}}.
\end{eqnarray*}
Similar to the argument in Lemma $4.2$, we can obtain the
results of Lemma $4.4$. $\Box$\\\\
\hspace*{0.6 cm} By $f^{+}\neq 0$, there exists at least one $u\in
W_{0}^{1,p}(\Omega)\backslash\{0\}$ such that
$\int_{\Omega}f|u|^{q}\,dx>0$. Let
$$\bar{t}_{\rm max}=\big(\frac{(2p-q)\lambda\int_{\Omega}f|u|^{q}\,dx}{pb\|u\|_{W^{1,p}}^{p}}\Big)^{\frac{1}{p-q}}.$$
Then we have the following result.\\\\
{\bf Lemma 4.5.} Let $r=2p$, $a<\frac{1}{\Lambda}$ and
$0<\lambda<\hat{\lambda}_{0}(a)$. Then for each $u\in
W_{0}^{1,p}(\Omega)$ with $\int_{\Omega}f|u|^{q}\,dx>0$, there is a
unique $0<t^{+}<\bar{t}_{\rm \max}$ such that $t^{+}u\in
\mathcal{N}_{\lambda,M}^{+}$ and
$$\mathcal{J}_{\lambda,M}(t^{+}u)=\inf_{0<t<\bar{t}_{\rm
\max}}\,\mathcal{J}_{\lambda,M}(tu)$$\\
{\bf Proof.} Fix $u\in W_{0}^{1,p}(\Omega)$ with
$\int_{\Omega}(f|u|^{q})\,dx>0$. Let
\begin{eqnarray*}
\bar{h}(t)=bt^{-p}\|u\|_{W^{1,p}}^{p}-t^{q-2p}\int_{\Omega}f|u|^{q}\,dx\,\text
{for} t>0.
\end{eqnarray*}
Clearly, $\bar{h}(t)\to -\infty$ as $t\to 0^{+}$ and $\bar{h}(t)\to
0$ as $t\to \infty$. Since
\begin{eqnarray*}
\bar{h}'(t)=-pbt^{-p-1}\|u\|_{W^{1,p}}^{p}+(2p-q)t^{q-2p-1}\int_{\Omega}f|u|^{q}\,dx,
\end{eqnarray*}
$\bar{h}'(t)=0$ at $t=\bar{t}_{\max}$,$\bar{h}'(t)>0$ for $t\in
[0,\bar{t}_{\max})$ and $\bar{h}'(t)<0$ for $t\in
(\bar{t}_{\max},\infty)$. Then $\bar{h}(t)$ achieves its maximum at
$\bar{t}_{\max}$, increasing for $t\in (0,\bar{t}_{\max})$ and
decreasing for $t\in (\bar{t}_{\max},\infty)$. Furthermore,
\begin{eqnarray*}
\bar{h}(\bar{t}_{\max})&=&\|u\|_{W^{1,p}}^{2p}b(\frac{p-q}{2p-q})(\frac{pb\|u\|_{W^{1,p}}^{q}}{(2p-q)\int_{\Omega}f|u|^{q}\,dx})^{\frac{p}{p-q}}\\
&\geq&=\|u\|_{W^{1,p}}^{2p}b(\frac{p-q}{2p-q})(\frac{pb\|u\|_{W^{1,p}}^{q}}{(2p-q)\|f^{+}\|_{\infty}})^{\frac{p}{p-q}}\\&\geq&(\frac{1-\Lambda
a}{\Lambda})\|u\|_{W^{1,p}}^{2p}.
\end{eqnarray*}
Since
$$\int_{\Omega}g|u|^{2p}\,dx-a\|u\|_{W^{1,p}}^{2p}\leq\,(\frac{1-\Lambda a}{\Lambda})\|u\|_{W^{1,p}}^{2p}\leq\,\bar{h}(\bar{t}_{max}),$$
and $\bar{h}(t)\to -\infty$ as $t\to 0^{+}$, we can conclude that
there is a unique $t^{+}<\bar{t}_{\max}$ such that
$\bar{h}(t^{+})=\int_{\Omega}g|u|^{2p}\,dx-a\|u\|_{W^{1,p}}^{2p}$
and $\bar{h}'(t^{+})>0$. The results of Lemma $4.5$ can be
obtained by repeating the argument of Lemma $4.2$.$\hspace* {.2cm}\Box$\\

To prove Theorem $4.1(iii)$, we require the
following lemma:\\\\
{\bf Lemma 4.6.} Suppose that $r=p^{2}$ and
$a\geq\frac{1}{\Lambda}$. Then for each $u\in W_{0}^{1,p}(\Omega)$
with $\int_{\Omega}f|u|^{q}\,dx>0$, there is unique
$0<t^{+}<\bar{t}_{\max}$ such that $t^{+}u\in
\mathcal{N}_{\lambda,M}^{+}$ and
$$\mathcal{J}_{\lambda,M}(t^{+}u)=\inf_{t\geq
0}\,\mathcal{J}_{\lambda,M}(tu).$$\\\\
{\bf Proof.} Similar to the argument in Lemma $4.5$, we can obtain
the results of lemma $4.6$.$\hspace* {.2cm}\Box$\\

To prove Theorem $4.1(iv)$, we use the
following lemma:\\\\
{\bf Lemma 4.7.} Suppose that $r<p^{2}$ and $a>\hat{A_{0}}$. Then
for each $u\in W_{0}^{1,p}(\Omega)$ with
$\int_{\Omega}f|u|^{q}\,dx>0$ and $\lambda>0$, there is a unique
$t_{\lambda}>0$ such that $t_{\lambda}u\in
\mathcal{N}_{\lambda,M}^{+}$ and
$$\mathcal{J}_{\lambda,M}(t_{\lambda}u)=\inf_{t\geq0}\,\mathcal{J}_{\lambda,M}(tu).$$\\\\
{\bf Proof.} Fix $u\in W_{0}^{1,p}(\Omega)$ with
$\int_{\Omega}f|u|^{q}\,dx>0$. Let
\begin{eqnarray*}
\tilde{h}(t)=at^{p^{2}-q}\|u\|_{W^{1,p}}^{p^{2}}+bt^{p-q}\|u\|_{W^{1,p}}^{p}-t^{r-q}\int_{\Omega}g|u|^{r}\,dx,\,\,\text
{for}\,\,t>0.
\end{eqnarray*}
Then by $q<p$ and $r<p^{2}$, we have$\tilde{h}(0)=0$ and
$\tilde{h}(t)\to \infty$ as$t\to \infty$. Since $a>\hat{A_{0}}$,
we have
\begin{align*}
\tilde{h}'(t)&=a(p^{2}-q)t^{p^{2}-q-1}\|u\|_{W^{1,p}}^{p^{2}}+b(p-q)t^{p-q-1}\|u\|_{W^{1,p}}^{p}-(r-q)t^{r-q-1}\int_{\Omega}g|u|^{r}\,dx\\
&\geq
t^{p-q-1}\|u\|_{W^{1,p}}^{p}(a(p^{2}-q)t^{p^{2}-p}+b(p-q)-(r-q)t^{r-q}\int_{\Omega}g|u|^{r}\,dx)\\&>0,\,\text
{for all}\,\,  t>0,
\end{align*}
and so $\tilde{h}(t)$ increases for $t\in[0,\infty)$. Moreover, for
each $\lambda>0$, there is a unique $t_{\lambda}>0$ such that
$$\tilde{h}(t_{\lambda})=\lambda\int_{\Omega}f|u|^{q}\,dx.$$
Again, the proof is completed by repeating the argument of Lemma
$4.2.$$\Box$
\section{ Proofs of Theorems 2.1, 2.2}
We write $\mathcal{N}_{\lambda,M}=\mathcal{N}_{\lambda,M}^{+}\cup
\mathcal{N}_{\lambda,M}^{-}$ and define
\begin{eqnarray*}
\alpha_{\lambda}^{+}=\inf_{u\in
\mathcal{N}_{\lambda,M}^{+}}\mathcal{J}_{\lambda,M}(u);
\,\,\,\alpha_{\lambda}^{-}=\inf_{u\in
\mathcal{N}_{\lambda,M}^{-}}\mathcal{J}_{\lambda,M}(u).
\end{eqnarray*}
Then we have the following result.\\\\
{\bf Theorem 5.1.} Suppose that $r>p^{2}$ and
$0<\lambda<\lambda_{0}(a)$. Then we have
\begin{itemize}
\item{\bf (i)} $\alpha_{\lambda}^{+}<0$;
\item{\bf (ii)} $\alpha_{\lambda}^{-}>c_{0}$ for some $c_{0}>0$.
\end{itemize}
In particular $\alpha_{\lambda}^{+}=\inf_{u\in
\mathcal{N}_{\lambda,M}}\mathcal{J}_{\lambda,M}(u)$.\\\\
{\bf Proof.} $(i)$ Let $u\in \mathcal{ N}_{\lambda,M}^{+}$. Since
\begin{eqnarray*}
\lambda(r-q)\int_{\Omega}f|u|^{q}\,dx&>&a(r-p^{2})\|u\|_{W^{1,p}}^{p^{2}}+b(r-p)\|u\|_{W^{1,p}}^{p}\\&\geq&
b(r-p)\|u\|_{W^{1,p}}^{p},
\end{eqnarray*}
then
\begin{eqnarray*}
\mathcal{J}_{\lambda,M}(u)&=&\frac{1}{p}\hat{M}(\|u\|_{W^{1,p}}^{p})-\frac{\lambda}{q}\int_{\Omega}f|u|^{q}\,dx-\frac{1}{r}\int_{\Omega}g|u|^{r}\,dx\\
&=&\frac{1}{p}\hat{M}(\|u\|_{W^{1,p}}^{p})-\frac{1}{p}M(\|u\|_{W^{1,p}}^{p})\,\|u\|_{W^{1,p}}^{p}\,-\lambda(\frac{r-q}{rq})\int_{\Omega}f|u|^{q}\,dx\\
&\leq&\frac{b(r-p)}{pr}\,\|u\|_{W^{1,p}}^{p}-\lambda(\frac{r-q}{rq})\int_{\Omega}f|u|^{q}\,dx\\&\leq&-\frac{b(r-p)(p-q)}{rpq}\,\|u\|_{W^{1,p}}^{p}<0.
\end{eqnarray*}
Thus, $\alpha_{\lambda}^{+}<0$.\\\\
$(ii)$ Let $u\in \mathcal{N}_{\lambda,M}^{-}$. We divide the proof
into the following two
cases.\\
Case (A): $r>p^{2}$ and $\lambda_{0}(a)=\frac{q\lambda_{2}}{p}$. By
$(4)$ and the Sobolev inequality,
\begin{eqnarray*}
b(p-q)\|u\|_{W^{1,p}}^{p}&\leq&a(p^{2}-q)\,\|u\|_{W^{1,p}}^{p^{2}}+b(p-q)\|u\|_{W^{1,p}}^{p}\\
&<&(r-q)S_{r}^{\frac{-r}{p}}\|g^{+}\|_{\infty}\,\|u\|_{W^{1,p}}^{r},
\end{eqnarray*}
this implies
$$\|u\|_{W^{1,p}}>(\frac{b(p-q)S_{r}^{\frac{r}{p}}}{(r-q)\|g^{+}\|_{\infty}})^{\frac{1}{r-p}}\,\,\,
\text {for all}  u\in \mathcal{N}_{\lambda,M}^{-}.$$ Subsequently,
\begin{eqnarray*}
\mathcal{J}_{\lambda }(
u)&\geq&\frac{\|u\|_{W^{1,p}}^{p}}{rp}\,(\frac{a(r-p^{2})}{p}\,\|u\|_{W^{1,p}}^{p^{2}-p}+b(r-p))-\lambda(\frac{r-q}{rq})\|f^{+}\|_{\infty}S_{q}^{-\frac{q}{p}}\|u\|_{W^{1,p}}^{q}\\
&\geq&\|u\|_{W^{1,p}}^{q}(\frac{b(r-p)}{rp}\|u\|_{W^{1,p}}^{p-q}-\lambda(\frac{r-q}{rq})\|f^{+}\|_{\infty}S_{q}^{-\frac{q}{p}})\\
&>&(\frac{b(p-q)S_{r}^{\frac{r}{p}}}{(r-p)\|g^{+}\|_{\infty}})^{\frac{q}{r-p}}\,(\frac{b(r-p)}{rp}(\frac{b(p-q)S_{r}^{\frac{r}{p}}}{(r-p)\|g^{+}\|_{\infty}})^{\frac{p-q}{r-p}}\,-\frac{\lambda(r-q)\|f^{+}\|_{\infty}}{rq
S_{q}^{\frac{q}{p}}}).
\end{eqnarray*}
Thus, if $\lambda<\frac{q}{p}\lambda_{2}$, then
$\alpha_{\lambda}^{-}>c_{0}$ for $c_{0}>0$.\\
Case (B): $r>p^{2}$ and
$\lambda_{0}(a)=\frac{q\lambda_{1}(a)}{p^{\frac{2p-1}{p}}}$. By
$(3)$ and the sobolev inequality,
\begin{eqnarray*}
p\sqrt[p]{ab^{p-1}(p^{2}-q)(p-q)^{p-1}}\|u\|_{W^{1,p}}^{2p-1}&\leq&a(p^{2}-q)\|u\|_{W^{1,p}}^{p^{2}}+b(p-q)(p-1)\|u\|_{W^{1,p}}^{p}
\\&<&(r-q)\|g^{+}\|_{\infty}S_{r}^{\frac{-r}{p}}\|u\|_{W^{1,p}}^{r},
\end{eqnarray*}
this implies
\begin{eqnarray*}
\|u\|_{W^{1,p}}>(\frac{pS_{r}^{\frac{r}{p}}\sqrt[p]{ab^{p-1}(p^{2}-q)(p-q)^{p-1}}}{(r-q)\|g^{+}\|_{\infty}})^{\frac{1}{(r-2p+1)}}\,\,\text
{for all}\,\,  u\in \mathcal{N}_{\lambda,M}^{-}.
\end{eqnarray*}
Repeating the argument of part $(A)$, we conclude that if
$\lambda<\frac{q\lambda_{1}(a)}{p^{\frac{2p-1}{p}}}$, then
$\alpha_{\lambda}^{-}>c_{0}$ for some $c_{0}>0$. this completes
the proof.$\hspace* {.2cm}\Box$\\\\
Now, we proceed to the proof of Theorem $2.1$. By Lemma $3.2$ and
the Ekeland variational principle \cite {ie}, there exist a
minimizing sequence $\{u_{n}^{\pm}\}$ for $\mathcal{J}_{\lambda ,M}$
on $\mathcal{N}_{\lambda,M}^{\pm}$ such that
\begin{eqnarray*}
\mathcal{J}_{\lambda
,M}(u_{n}^{\pm})=\alpha_{\lambda,M}^{\pm}+o(1)\,\,\text {and}\,\,
\mathcal{J}'_{\lambda,M}(u_{n}^{\pm})=o(1)\,\,\text {in}\,\,
W^{-1,p'}(\Omega).
\end{eqnarray*}
It follows, by Lemma $3.1$, that there exists a subsequence
$\{u_{n}^{\pm}\}$ and $u_{0}^{\pm}\in W_{0}^{1,p}(\Omega)$ are
solutions of Eq. $(1)$ such that $u_{n}^{\pm}\to u_{0}^{\pm}$
strongly in $W_{0}^{1,p}(\Omega)$ and so $u_{0}^{\pm}\in
\mathcal{N}_{\lambda,M}^{\pm}$ and $\mathcal{J}_{\lambda
,M}(u_{0}^{\pm})=\hat{\alpha}_{\lambda,M}^{\pm}$. Since
$\mathcal{J}_{\lambda ,M}(u_{0}^{\pm})=\mathcal{J}_{\lambda
,M}(|u_{0}^{\pm}|)$ and $|u_{0}^{\pm}|\in
\mathcal{N}_{\lambda,M}^{\pm}$, by Lemma $2.1$, we may assume that
$u_{0}^{\pm}$ are positive solutions of Eq. $(1)$. Moreover,
$\mathcal{N}_{\lambda,M}^{+}\cap \mathcal{N}_{\lambda,M}^{-}=\O$,
this indicates that $u_{0}^{+}$ and $u_{0}^{-}$ are two distinct
solutions.This completes the proof.\\
To prove Theorem $2.2$, we need the following.\\
By Theorem $4.3$, we write
$\mathcal{N}_{\lambda,M}=\mathcal{N}_{\lambda,M}^{+}\cup
\mathcal{N}_{\lambda,M}^{-}$ and define
\begin{eqnarray*}
\hat{\alpha}_{\lambda,M}^{+}=\inf_{u\in
\mathcal{N}_{\lambda,M}^{+}}\mathcal{J}_{\lambda,M}(u);
\,\,\,\hat{\alpha}_{\lambda,M}^{-}=\inf_{u\in
\mathcal{N}_{\lambda,M}^{-}}\mathcal{J}_{\lambda,M}(u).
\end{eqnarray*}
Then we have the following result.\\\\
 {\bf Theorem 5.2.} Suppose
that $r=p^{2}$, $a<\frac{1}{\Lambda}$
and $0<\lambda<\frac{1}{p}\hat{\lambda}_{0}(a)$. Then we have\\
\begin{itemize}
\item{\bf (i)} $\hat{\alpha}_{\lambda,M}^{+}<0$;
\item{\bf (ii)} $\hat{\alpha}_{\lambda,M}^{-}>c_{0}$ for some
$c_{0}>0$.
\end{itemize}
In particular, $\hat{\alpha}_{\lambda,M}^{+}=\inf_{u\in
\mathcal{N}_{\lambda,M}}\mathcal{J}_{\lambda,M}(u)$.\\\\
{\bf Proof.} $(i)$ Let $u\in \mathcal{N}_{\lambda,M}^{+}$. Since
$$\lambda(p^{2}-q)\int_{\Omega}f|u|^{q}\,dx>b(p^{2}-p)\|u\|_{W^{1,p}}^{p},$$
then
\begin{eqnarray*}
\mathcal{J}_{\lambda,M}(u)&=&\frac{1}{p}\hat{M}(\|u\|_{W^{1,p}}^{p})-\frac{\lambda}{q}\int_{\Omega}f|u|^{q}\,dx-\frac{1}{r}\int_{\Omega}g|u|^{r}\,dx
\\&=&\frac{1}{p}\hat{M}(\|u\|_{W^{1,p}}^{p})-\frac{1}{r}M(\|u\|_{W^{1,p}}^{p})\,\|u\|_{W^{1,p}}^{p}\,-\lambda(\frac{r-q}{rq})\int_{\Omega}f|u|^{q}\,dx
\\&=&\frac{b(p-1)}{p^{2}}\,\|u\|_{W^{1,p}}^{p}-\lambda(\frac{p^{2}-q}{p^{2}q})\int_{\Omega}f|u|^{q}\,dx
\\&<&-\frac{b(p-1)(p-q)}{p^{2}q}\,\|u\|_{W^{1,p}}^{p}<0.
\end{eqnarray*}
Thus, $\hat{\alpha}_{\lambda,M}^{+}<0$.\\
\hspace{0.6 cm}$(ii)$ Let $u\in \mathcal{N}_{\lambda,M}^{-}$. By
$(4)$,
\begin{eqnarray*}
b(p-q)\|u\|_{W^{1,p}}^{p}&<&(p^{2}-q)\int_{\Omega}g|u|^{r}\,dx-a(p^{2}-q)\|u\|_{W^{1,p}}^{p^{2}}\\&\leq&\frac{(p^{2}-q)(1-a\Lambda)}{\Lambda}\|u\|_{W^{1,p}}^{p^{2}},
\end{eqnarray*}
which implies that
\begin{equation}
\|u\|_{W^{1,p}}>(\frac{b\Lambda\,(p-q)}{(1-a\Lambda)(p^{2}-q)})^{\frac{1}{p^{2}-p}}\,\,\text
{for all}  u\in \mathcal{N}_{\lambda,M}^{-}.
\end{equation}
Subsequently,
\begin{equation}
\begin{aligned}
\mathcal{J}_{\lambda,M}(u)&=\frac{1}{p}
\hat{M}(\|u\|_{W^{1,p}}^{p})-\frac{\lambda}{q}\int_{\Omega}f|u|^{q}\,dx-\frac{1}{r}\int_{\Omega}g|u|^{r}\,dx\\
&\geq
\frac{b(p-1)}{p^{2}}\,\|u\|_{W^{1,p}}^{p}-\lambda(\frac{p^{2}-q}{p^{2}q})
\|f^{+}\|_{\infty}S_{q}^{\frac{-q}{p}}\|u\|_{W^{1,p}}^{q}\\&\geq
\|u\|_{W^{1,p}}^{q}(\frac{b(p-1)}{p^{2}}\,\|u\|_{W^{1,p}}^{p-q}-\lambda(\frac{p^{2}-q}{p^{2}q})\|f^{+}\|_{\infty}S_{q}^{\frac{-q}{p}})
\\&> \big(\frac{b\Lambda (p-q)}{(1-a\Lambda)(p^{2}-q)}\Big)^{\frac{q}{p^{2}-p}}\big(\frac{b(p-1)}
{p^{2}}\big(\frac{b\Lambda (p-q)}{(1-a\Lambda)(p^{2}-q)})\Big)^{\frac{p-q}{p^{2}-p}}-\lambda(\frac{p^{2}-q}{p^{2}q})\|f^{+}\|_{\infty}S_{q}^{\frac{-q}{p}}\Big).
\end{aligned}
\end{equation}
Thus, if $\lambda<\frac{1}{p}\hat{\lambda}_{0}(a)$, then
$\alpha_{\lambda}^{-}>c_{0}$ for some $c_{0}>0$. This completes
the proof.$\Box$\\

Now, we proceed to the proof of Theorem $2.2.$$(i)$ By Theorem
$4.1$$(iii)$, we define
$$\hat{\alpha}_{\lambda,M}=\inf_{u\in
\mathcal{N}_{\lambda,M}^{+}}\mathcal{J}_{\lambda,M}(u).$$ Similar to
the argument in Theorem $5.3$, we can conclude
$\hat{\alpha}_{\lambda,M}<0$. Moreover, by Lemma $3.2(i)$ and the
Ekeland variational principle \cite {ie}, there exist a minimizing
sequence $\{u_{n}\}$ for $\mathcal{J}_{\lambda ,M}$ on
$\mathcal{N}_{\lambda,M}^{+}$ such that
\begin{eqnarray*}
\mathcal{J}_{\lambda ,M}(u_{n})=\alpha_{\lambda,M}+o(1)\text {and}
\mathcal{J}'_{\lambda,M}(u_{n})=o(1)\text {in} W^{-1,p'}(\Omega).
\end{eqnarray*}
It follows, by Lemma $3.1$, that there exists a subsequence
$\{u_{n}\}$ and $u_{0}\in W_{0}^{1,p}(\Omega)$ is a solution of Eq.
$(2)$ such that $u_{n}\to u_{0}$ strongly in $W_{0}^{1,p}(\Omega)$
and so $u_{0}\in \mathcal{N}_{\lambda,M}^{+}$ and
$\mathcal{J}_{\lambda ,M}(u_{0})=\hat{\alpha}_{\lambda}$. Since
$\mathcal{J}_{\lambda ,M}(u_{0})=\mathcal{J}_{\lambda ,M}(|u_{0}|)$
and $|u_{0}|\in \mathcal{N}_{\lambda,M}^{+}$, by lemma $2.1$, we may
assume that $u_{0}$ is a positive solution of
Eq. $(1)$.\\
$(ii)$ Similar to the argument in Theorem $2.1$, Eq. $(1)$ thus has
two positive solutions $u_{\lambda,M}^{+}\in
\mathcal{N}_{\lambda,M}^{+}$ and $u_{\lambda,M}^{-}\in
\mathcal{N}_{\lambda,M}^{-}$. Moreover, by $(18)$ and $(19)$,
$$\|u_{\lambda,M}^{-}\|_{W^{1,p}}\to \infty\,\,\,\text {as}\,\,  a\to
\frac{1^{-}}{\Lambda},$$ and
$$\lim_{a\to \frac{1^{-}}{\Lambda}}\,\inf_{u\in
\mathcal{N}_{\lambda,M}^{-}}\,\mathcal{J}_{\lambda ,M}(u)=\infty.$$
This completes the proof.
\section{ Proof of Theorem 2.3}
\hspace{0.6 cm} First, we consider the following truncated equation
$(1)$,
\begin{equation}
\left\{\begin{array}{ll} -M_{k^{p-1}}\,\Big(\int_{\Omega }|\nabla
u|^{p}\,dx\Big) \Delta_{p}u=\lambda f(x)|u|^{q-2}u+g(x)|u|^{r-2}u ,&
\in \Omega,\\ u=0 , & \in
\partial \Omega.
\end{array}\right.
\end{equation}
where ${k^{p-1}}\in (\frac{b(r-p)}{ra},\frac{b(r-p)}{pa})$ and
\begin{eqnarray*}
M_{k^{p-1}}(s)=
\begin{cases}
M(s), & if s\leq k^{p-1},\\  M(k^{p-1}) &  if s>k^{p-1},
\end{cases}
\end{eqnarray*}
is a truncated function of $M(s)$. The positive solutions of
truncated equation $(20)$ are critical points of the functional
\begin{eqnarray*}
\mathcal{J}_{\lambda,M_{k^{p-1}}}(u)=\frac{1}{p}\hat{M}_k^{p-1}(\|u\|_{W^{1,p}}^{p})-\frac{\lambda}{q}\int_{\Omega}f|u|^{q}\,dx-\frac{1}{r}\int_{\Omega}g|u|^{r}\,dx,
\end{eqnarray*}
where
$\hat{M}_{k^{p-1}}(\|u\|_{W^{1,p}}^{p})=\int_{0}^{t}M_{k^{p-1}}\,ds$.
We have the following results.\\\\
{\bf Lemma 6.1.} The energy functional
$\mathcal{J}_{\lambda,M_{k^{p-1}}}$ is coercive and
bounded below on $\mathcal{N}_{\lambda,M_{k^{p-1}}}$.\\\\
{\bf Proof.} For $u\in \mathcal{N}_{\lambda,M_{k^{p-1}}}$, we have
$M_{k^{p-1}}\,(\|u\|_{W^{1,p}}^{p})=\lambda
\int_{\Omega}f|u|^{q}\,dx+\int_{\Omega}g|u|^{r}\,dx$. By the
sobolev inequality,
\begin{eqnarray*}
\mathcal{J}_{\lambda,M_{k^{p-1}}}(u)&=&\frac{1}{p}\hat{M}_{k^{p-1}}(\|u\|_{W^{1,p}}^{p})-\frac{\lambda}{q}\int_{\Omega}f|u|^{q}\,dx-\frac{1}{r}\int_{\Omega}g|u|^{r}\,dx
\\&=&\frac{1}{p}\hat{M}_{k^{p-1}}(\|u\|_{W^{1,p}}^{p})-\frac{1}{r}\,M_{k^{p-1}}\,(\|u\|_{W^{1,p}}^{p})\,\|u\|_{W^{1,p}}^{p}\,-\lambda(\frac{r-q}{rq})\int_{\Omega}f|u|^{q}\,dx
\end{eqnarray*}
\begin{equation}
\hspace*{-1
cm}\geq\,(\frac{b}{p}-\frac{M_{k^{p-1}}}{r})\|u\|_{W^{1,p}}^{p}-\lambda(\frac{r-q}{rq})\|f^{+}\|_{\infty}S_{q}^{\frac{-q}{p}}\|u\|_{W^{1,p}}^{q},
\end{equation}
and since $k^{p-1}<\frac{b(r-p)}{pa}$, this gives
$\frac{b}{p}\,-\frac{M(k^{p-1})}{r}\,>0$. Thus,
$\mathcal{J}_{\lambda,k^{p-1}}$ is coercive and bounded below
on $\mathcal{N}_{\lambda,M_{k^{p-1}}}$.$\hspace* {.2cm}\Box$\\

By $(4)$ and $(5)$, if $u\in \mathcal{N}_{\lambda,M_{k^{p-1}}}$ with
$\|u\|_{W^{1,p}}^{p}\leq\,k^{p-1}$,
then\\
\begin{eqnarray*}
I''_{u,M_{k^{p-1}}}\,(1)=a(p^{2}-1)\,\|u\|_{W^{1,p}}^{p^{2}}+(p-1)b\|u\|_{W^{1,p}}^{p}-\lambda\,(q-1)\int_{\Omega}f|u|^{q}\,dx-(r-1)\int_{\Omega}g|u|^{r}\,dx
\end{eqnarray*}
\begin{equation}
\hspace*{-1.7
cm}=[a(p^{2}-q)\|u\|_{W^{1,p}}^{p^{2}-p}+b(p-q)]\,\|u\|_{W^{1,p}}^{p}-(r-q)\int_{\Omega}g|u|^{r}\,dx
\end{equation}
\begin{equation}
\hspace*{-2.5
cm}=[a(p^{2}-r)\|u\|_{W^{1,p}}^{p^{2}-p}+b(p-r)]\,\|u\|_{W^{1,p}}^{p}+\lambda\,\int_{\Omega}f|u|^{q}\,dx,
\end{equation}
and if $u\in \mathcal{N}_{\lambda,M_{k^{p-1}}}$ with
$\|u\|_{W^{1,p}}^{p}>k^{p-1}$, then
\begin{align*}
I''_{u,M_{k^{p-1}}}\,(1)=(p-1)\,M(k^{p-1})\,\|u\|_{W^{1,p}}^{p}-\lambda\,(q-1)\int_{\Omega}f|u|^{q}\,dx-(r-1)\int_{\Omega}g|u|^{r}\,dx
\end{align*}
\begin{equation}
\hspace*{-2
cm}=(p-q)\,M(k^{p-1})\,\|u\|_{W^{1,p}}^{p}-(r-q)\int_{\Omega}g|u|^{r}\,dx
\end{equation}
\begin{equation}
\hspace*{-1
cm}=-(r-p)\,M(k^{p-1})\,\|u\|_{W^{1,p}}^{p}\,+\lambda\,(r-q)\,\int_{\Omega}f|u|^{q}\,dx.
\end{equation}
Furthermore, if $u\in \mathcal{N}_{\lambda,M_{k^{p-1}}}^{0}$, by
$(22)$-$(25)$ and the Sobolev inequality, then
\begin{equation}
\tilde{C_{1}}\|u\|_{W^{1,p}}^{p}\leq\,(r-q)\int_{\Omega}g|u|^{r}\,dx\leq\,(r-q)S_{r}^{\frac{-r}{p}}\|g^{+}\|_{\infty}\|u\|_{W^{1,p}}^{r},
\end{equation}
where $\tilde{C}_{1}=(p-q)\,\min\{b,M(k^{p-1})\}$, and
\begin{equation}
\begin{aligned}
((r-p)b-a(p^{2}-r)k^{p-1})\|u\|_{W^{1,p}}^{p}\,&\leq
[a(r-p^{2})\|u\|_{W^{1,p}}^{p^{2}-p}+(r-p)]\,\|u\|_{W^{1,p}}^{p}\\
&\leq\lambda\,(r-q)\|f^{+}\|_{\infty}S_{q}^{\frac{-q}{p}}\,\|u\|_{W^{1,p}}^{q},\,\,\text
 {if}\,\, \|u\|_{W^{1,p}}^{p}\leq\,k^{p-1}.
\end{aligned}
\end{equation}
Note that $b(r-p)-a(p^{2}-r)k^{p-1}>0$ since
$k^{p-1}<\frac{b(r-p)}{pa}<\frac{b(r-p)}{a(p^{2}-r)}$, and so with
$(27)$,
\begin{equation}(b(r-p)-a(p^{2}-r)k^{p-1})\|u\|_{W^{1,p}}^{p}\leq\lambda\,(r-q)\|f^{+}\|_{\infty}\,S_{q}^{\frac{-q}{p}}\|u\|_{W^{1,p}}^{q}\,\,\text
 {if}  \|u\|_{W^{1,p}}^{p}\leq\,k^{p-1}.
\end{equation}
Moreover, by $(25)$,
\begin{equation}
\begin{aligned}
M(k^{p-1})(r-p)\|u\|_{W^{1,p}}^{p}&=\lambda\,(r-q)\int_{\Omega}f|u|^{q}\,dx\\
&\leq\lambda\,(r-q)\|f^{+}\|_{\infty}\,S_{q}^{\frac{-q}{p}}\,\|u\|_{W^{1,p}}^{q},\,\,\text
{if}\,\, \|u\|_{W^{1,p}}^{p}>k^{p-1}.
\end{aligned}
\end{equation}
It follows that, by $(28)$ and $(29)$,
\begin{equation}
\tilde{C_{2}}\|u\|_{W^{1,p}}^{p}\leq
\lambda(r-q)\|f^{+}\|_{\infty}S_{q}^{\frac{-q}{p}}\|u\|_{W^{1,p}}^{q},
\end{equation}
where
$\tilde{C_{2}}=\min\{M(k^{p-1})\,(r-p),[b(r-p)-a(p^{2}-p)k^{p-1}]\}$.
Hence, by $(26)$ and $(30)$,
\begin{eqnarray*}
(\frac{S_{r}^{\frac{r}{p}}\tilde{C}_{1}}{(r-q)\|g^{+}\|_{\infty}})^{\frac{1}{r-p}}\leq\|u\|_{W^{1,p}}\leq\,(\frac{\lambda\,(r-q)\|f^{+}\|_{\infty}}{S_{q}^{\frac{q}{p}}\tilde{C_{2}}})^{\frac{1}{p-q}},
\end{eqnarray*}
for all $u\in \mathcal{N}_{\lambda,M_{k^{p-1}}}^{0}$. Thus, if
the submanifold $\mathcal{N}_{\lambda,M_{k^{p-1}}}^{0}$ is
nonempty, then the inequality
\begin{eqnarray*}
\lambda\geq\,\tilde{C_{2}}\,\tilde{C_{3}}\,\,\text {where}\,\,
  \tilde{C}_{3}=(\frac{S_{r}^{\frac{r}{p}}
  \tilde{C}_{1}}{(r-q)\|g^{+}\|_{\infty}})^{\frac{p-q}{r-p}}\frac{S_{q}^{\frac{q}{p}}}{(r-q)\|f^{+}\|_{\infty}},
\end{eqnarray*}
must be hold. Subsequently, we have the following result.\\\\
{\bf Lemma 6.2.} If $0<\lambda<\tilde{C}_{2}\,\tilde{C}_{3}$,then
the submanifold $\mathcal{N}_{\lambda,M_{k^{p-1}}}^{0}=\emptyset$.\\

By Lemma $(6.2)$, we write
$\mathcal{N}_{\lambda,M_{k^{p-1}}}=\mathcal{N}_{\lambda,M_{k^{p-1}}}^{+}\cup\,\mathcal{N}_{\lambda,M_{k^{p-1}}}^{-}$.
Using a similar argument to that of Lemma $4.3$, it can be deduced
that $\mathcal{N}_{\lambda,M_{k^{p-1}}}^{\pm}\,\neq\emptyset$.
Define
\begin{eqnarray*}
\alpha_{\lambda,M_{k^{p-1}}}^{+}=\inf_{u\in
\mathcal{N}_{\lambda,M_{k^{p-1}}}^{+}}\,\mathcal{J}_{\lambda,M_{k^{p-1}}}(u);\,\,\,
\alpha_{\lambda,M_{k^{p-1}}}^{-}=\inf_{u\in
\mathcal{N}_{\lambda,M_{k^{p-1}}}^{-}}\,\mathcal{J}_{\lambda,M_{k^{p-1}}}(u),
\end{eqnarray*}
then we have the following result.\\\\
{\bf Theorem 6.3.} We have
\begin{itemize}
\item[(i)] $\alpha_{\lambda,M_{k^{p-1}}}^{+}<0$ for all $\lambda\in (0,\tilde{C}_{2}\tilde{C}_{3})$;
\item[(ii)] if $0<\lambda<
\tilde{C}_{3}\tilde{C}_{4}$, then
$\alpha_{\lambda,k^{p-1}}^{-}>c_{0}$ for some $c_{0}>0$, where
$\tilde{C_{4}}=\frac{q(rb-pM(k^{p-1}))}{p}$.
\end{itemize}
In particular, for each
$0<\lambda<\tilde{C}_{3}\,\min\{\tilde{C}_{2},\tilde{C}_{4}\}$, we
have $$\alpha_{\lambda,k^{p-1}}^{+}=\inf_{u\in
\mathcal{N}_{\lambda,M_{k^{p-1}}}}\,\mathcal{J}_{\lambda,M_{k^{p-1}}}(u)$$.\\
{\bf Proof}. $(i)$ Let $u\in
\mathcal{N}_{\lambda,M_{k^{p-1}}}^{+}$. We
divide the proof into the following three cases.\\
Case (A): $\|u\|_{W^{1,p}}^{p}\leq\,k^{p-1}$. By $(23)$,
\begin{equation}
\begin{aligned}
0<(b(r-p)-a(p^{2}-r)k^{p-1})\|u\|_{W^{1,p}}^{p}&\leq\,[a\,(r-p^{2})\,\|u\|_{W^{1,p}}^{p^{2}-p}\,+b(r-p)]\,\|u\|_{W^{1,p}}^{p}\\
&<\lambda\,(r-q)\,\int_{\Omega}f|u|^{q}\,dx.
\end{aligned}
\end{equation}
Since $b\,(r-p)\,-a\,(p^{2}-r)\,k^{p-1}>0$, it follows that
\begin{align*}
\mathcal{J}_{\lambda,M_{k^{p-1}}}(u)&=\frac{1}{p}\hat{M}_{k^{p-1}}(\|u\|_{W^{1,p}}^{p})-\frac{1}{r}\,M_{k^{p-1}}\,(\|u\|_{W^{1,p}}^{p})\,\|u\|_{W^{1,p}}^{p}\,-\lambda(\frac{r-q}{rq})\int_{\Omega}f|u|^{q}\,dx
\\&=\frac{a(r-p^{2})}{rp^{2}}\|u\|_{W^{1,p}}^{p^{2}}\,+\frac{b(r-p)}{rp}\|u\|_{W^{1,p}}^{p}\lambda(\frac{r-q}{rq})\int_{\Omega}f|u|^{q}\,dx\\&<\frac{\|u\|_{W^{1,p}}^{p}}{rpq}[\frac{a(p^{2}-q)(p^{2}-r)}{p}\|u\|_{W^{1,p}}^{p^{2}-p}-b(r-p)(p-q)]
\\&\leq\,-\frac{k^{p-1}}{prq}\,[b(p-q)(r-p)-\frac{a(p^{2}-q)(p^{2}-r)}{p}k^{p-1}]\\&<-\frac{k^{p-1}}{prq}\,[b(r-p)-a(p^{2}-r)\,k^{p-1}]<0.
\end{align*}
Case (B): $\|u\|_{W^{1,p}}^{p}>k^{p-1}$. By $(25)$,
$$
M(k^{p-1})\,(r-p)\,\|u\|_{W^{1,p}}^{p}<\lambda\,(r-q)\int_{\Omega}f|u|^{q}\,dx.
$$
Moreover,
\begin{eqnarray*}
\hat{M}_{k^{p-1}}(t)&=&\int_{0}^{t}M_{k^{p-1}}(s)\,ds=\int_{0}^{k^{p-1}}M_{k^{p-1}}(s)\,ds+\int_{k^{p-1}}^{t}M_{k^{p-1}}(s)\,ds
\\&=&\int_{0}^{k^{p-1}}M(s)\,ds+\int_{k^{p-1}}^{t}M(s)\,ds\\&=&\hat{M}(k^{p-1})+M(k^{p-1})(t-k^{p-1})
\,\text {for} t>k^{p-1},
\end{eqnarray*}
and thus,
\begin{align*}
\mathcal{J}_{\lambda,M_{k^{p-1}}}(u)&=\frac{1}{p}\,\hat{M}_{k^{p-1}}(\|u\|_{W^{1,p}}^{p})-\frac{1}{r}\,M_{k^{p-1}}\,(\|u\|_{W^{1,p}}^{p})\,\|u\|_{W^{1,p}}^{p}\,-\lambda(\frac{r-q}{rq})\,\int_{\Omega}f|u|^{q}\,dx
\\&=\frac{1}{p}\,(\hat{M(k^{p-1})}\,-\,M\,(k^{p-1})\,k^{p-1})+\frac{M\,(k^{p-1})\,(r-p)}{rp}\,\|u\|_{W^{1,p}}^{p}\\&-\lambda(\frac{r-q}{rq})\,\int_{\Omega}f|u|^{q}\,dx
\\&<\frac{1}{p}\,(\frac{1}{p}\,a\,k^{p(p-1)}\,+\,b\,k^{p-1}\,-a\,k^{2(p-1)}\,-b\,k^{p-1})+\frac{\lambda(r-q)}{rp}\int_{\Omega}f|u|^{q}\,dx\\&-\lambda(\frac{r-q}{rq})\int_{\Omega}f|u|^{q}\,dx
\\&=-\frac{1}{p^{2}}\,a\,k^{p-1}\,(p\,k^{p-1}-k^{(p-1)^{2}})-\frac{\lambda\,(p-q)(r-q)}{rpq}\,\int_{\Omega}f|u|^{q}\,dx<0.
\end{align*}
Consequently, $\alpha_{\lambda,M_{k^{p-1}}}^{+}<0$.\\
$(ii)$ Let $u\in \mathcal{N}_{\lambda\,M(k^{p-1})}^{-}$. By
$(22)$,$(24)$ and the Sobolev inequality,
$$
\min\{b,M(k^{p-1})\}\,(p-q)\,\|u\|_{W^{1,p}}^{p}<(r-q)\,\int_{\Omega}g|u|^{r}\,dx\leq\,(r-q)\,\|g^{+}\|_{\infty}\,S_{r}^{\frac{-r}{p}}\|u\|_{W^{1,p}}^{r},
$$
this implies
\begin{equation}
\|u\|_{W^{1,p}}\,>S_{r}^{\frac{r}{p(r-p)}}\,\Big(\frac{\min\{b,M(k^{p-1})\}\,(p-q)}{\|g^{+}\|_{\infty}(r-q)\,}\Big)^{\frac{1}{(r-p)}}
\,\,\text  {for all}\,\,\,  u\in
\mathcal{N}_{\lambda\,M_{k^{p-1}}}^{-}.
\end{equation}
By $(20)$ from the proof of Lemma $6.3$,
\begin{align*}
\mathcal{J}_{\lambda,M_{k^{p-1}}}(u)&\geq
\|u\|_{W^{1,p}}^{q}\Big[\frac{br\,-p\,M(k^{p-1})}
{pr}\,\|u\|_{W^{1,p}}^{p-q}\,-\frac{\lambda\,(r-q)}{rqS_{q}^{\frac{q}{p}}}\,\|f^{+}\|_{\infty}\Big]
\\&>S_{r}^{\frac{rq}{p(r-p)}}\,\Big(\frac{\min\{b,M(k^{p-1})\}\,(p-q)}{\|g^{+}\|_{\infty}(r-q)}\Big)^{\frac{q}{(r-p)}}\\
&\Big(\frac{(br\,-p\,M(k^{p-1}))S_{r}^{\frac{r(p-q)}{p(r-p)}})}{pr}\,(\frac{\min\{b,M(k^{p-1})\}\,(p-q)}{\|g^{+}\|_{\infty}(r-q)})^{\frac{p-q}{(r-p)}}\,-\frac{\lambda\,(r-q)}{rqS_{q}^{\frac{q}{p}}}\,\|f^{+}\|_{\infty}\Big).
\end{align*}
Thus, if $\lambda\,<\tilde{C}_{3}\,\tilde{C}_{4}$, then
$\alpha_{\lambda,M_{k^{p-1}}}^{-}>c_{0}$. This completes the
proof.\hspace*
{0.4 cm}$\Box$\\

Now, we proceed to the proof of Theorem 2.3.(i) By Lemma 3.2 (ii)
and the Ekeland variational principle \cite {ie}, there exists a
minimizing sequence $\{u_{n}\}$ for $\mathcal{J}_{\lambda,M}$ on
$W_{0}^{1,p}(\Omega)$ such that
$$
\mathcal{J}_{\lambda\,M}\,(u_{n})=\beta_{\lambda}\,+\,o(1)\,\,\text
{and} \mathcal{J}'_{\lambda\,M}\,(u_{n})=\,o(1)\,\text {in}\,\,
W^{-1,p'},
$$
where $\beta_{\lambda}=\inf_{u\in
W_{0}^{1,p}(\Omega)}\mathcal{J}_{\lambda,M}(u)$. Clearly,
$\beta_{\lambda}<0$. Then by Lemma $3.1$, there exist a subsequence
$\{u_{n}\}$ and $u_{a,\lambda}\in W_{0}^{1,p}(\Omega)$ is a nonzero
solution of Eq. $(1)$ such that $u_{n}\to u_{0}$ strongly in
$W_{0}^{1,p}\,(\Omega)$ and
$\mathcal{J}_{\lambda,M}\,(u_{a,\lambda})=\mathcal{J}_{\lambda,M}\,(|u_{a,\lambda}|)$,
by Lemma 2.1 we may assume that $u_{a,\lambda}$ is a positive
solution of Eq. $(1)$.\\
\hspace{0.6 cm} $(ii)$ Let $\theta\,>0$ and take
$\lambda\,<\tilde{\lambda}_{0}=\min\{\theta,\tilde{C}_{3}\,\min\{\tilde{C_{2}},\tilde{C}_{4}\}\}$.
Then by Lemma $6.1$ and the Ekeland variational principle \cite
{ie}, there exist two minimizing sequences $\{u_{n}^{\pm}\}$ for
$\mathcal{J}_{\lambda,M_{k^{p-1}}}$ on
$\mathcal{N}_{\lambda,M_{k^{p-1}}}^{\pm}$ such that
$$
\mathcal{J}_{\lambda,M_{k^{p-1}}}\,(u_{n}^{\pm})=\alpha_{\lambda,k^{p-1}}^{\pm}\,+\,o(1)
\,\,\text {and}\,\,
 \mathcal{J}'_{\lambda,M_{k^{p-1}}}\,(u_{n}^{\pm})=\,o(1)\,\,\text
 {in}\,\,  W^{-1,p'}.
$$
Using a similar argument to that in Lemma $3.1$, there exist
subsequences $\{u_{n}^{\pm}\}$ and
${u_{\lambda,M_{k^{p-1}}}^{\pm}}\,\in W_{0}^{1,p}\,(\Omega)$ are
nonzero solutions of Eq. $(20)$ such that
${u_{n}^{\pm}}\rightharpoonup\,{u_{\lambda}^{\pm}}$ strongly in
$W_{0}^{1,p}\,(\Omega)$ and so $u_{\lambda,M_{k^{p-1}}}^{\pm}\,\in
\mathcal{N}_{\lambda,M_{k^{p-1}}}^{\pm}$ and
$\mathcal{J}_{\lambda,M_{k^{p-1}}}(u_{\lambda,M_{k^{p-1}}}^{\pm})=\alpha_{\lambda\,k^{p-1}}^{\pm}$.
Since
$\mathcal{J}_{\lambda\,M_{k^{p-1}}},(u_{{\lambda,M_{k^{p-1}}}^{\pm}})=\mathcal{J}_{\lambda,M_{k^{p-1}}}\,(|u_{{\lambda,M_{k^{p-1}}}^{\pm}}|)$
and $|u_{\lambda\,M_{k^{p-1}}}^{\pm}|\,\in
\mathcal{N}_{\lambda,M_{k^{p-1}}}^{\pm}$, by Lemma $2.1$, it can be
deduced that $u_{\lambda,M_{k^{p-1}}}^{\pm}$ are positive solutions
of Eq. $(20)$. Moreover,
$N_{\lambda,M_{k^{p-1}}}^{+}\,\cap\,\mathcal{N}_{\lambda,M_{k^{p-1}}}^{-}=\emptyset$,
this implies that $u_{\lambda,M_{k^{p-1}}}^{+}$ and
$u_{\lambda,M_{k^{p-1}}}^{-}$ are two distinct solutions. Now, we
claim that
$\|u_{\lambda,M_{k^{p-1}}}^{\pm}\|_{W^{1,p}}^{p}\,\leq\,k^{p-1}$; if
this is not the case, then by $(11)$ and
$k\in\,(\frac{b(r-p)}{ra},\frac{b(r-p)}{pa})$,
\begin{eqnarray*}
\frac{b(r-p)}{ar\,L(\theta)}&=&\frac{b(r-p)}{ra(\theta\,C_{*}^{q}\,\|f^{+}\|_{\infty}+C_{*}^{r}\,\|g^{+}\|_{\infty})\,|\Omega|}
\\&<&\frac{k^{p-1}}{(\lambda\,C_{*}^{q}\,\|f^{+}\|_{\infty}+C_{*}^{r}\,\|g^{+}\|_{\infty})\,|\Omega|}
\\&<&\max\{(\frac{b(2r-p)}{r})^{\frac{(p-r+q)}{(r-p)}},(\frac{br}{p})^{\frac{(p-r+q)}{(r-p)}},(\frac{br}{p})^{\frac{p}{(r-p)}}\}=A_{0},
\end{eqnarray*}
which implies $a>\frac{b(r-p)}{r\,A_{0}\,L(\theta)}$ a
contradiction. Thus,
$u_{{\lambda,M_{k^{p-1}}}}^{\pm}=u_{{\lambda,M}}^{\pm}\,\in
\mathcal{N}_{\lambda,M}^{\pm}$ and
$\mathcal{J}_{\lambda,M}\,(u_{{\lambda\,M_{k^{p-1}}}}^{\pm})=\mathcal{J}_{\lambda,M_{k^{p-1}}}\,(u_{{\lambda,M_{k^{p-1}}}}^{\pm})=\alpha_{\lambda,(k^{p-1})}^{\pm}$.
Moreover,$u_{\lambda,M}^{+}$ and $u_{\lambda,M}^{-}$ are positive
solutions of Eq. $(1)$.
\section{Proof of Theorem 2.4}
\hspace*{0.6 cm} First, we consider a modified version of Eq. $(1)$
as follows,
\begin{equation}
\left\{\begin{array}{ll} -M_{\hat{k}}\,\Big(\int_{\Omega }|\nabla
u|^{p}dx\Big) \Delta_{p}u=\lambda f(x)|u|^{q-2}u+g(x)|u|^{r-2}u ,&
\in \Omega,\\ u=0 , & \in
\partial \Omega.
\end{array}\right.
\end{equation}
where $\hat{k}=\frac{b(r-p)}{a(2p-r)}$ and \\
\begin{eqnarray*}
M_{\hat{k}}(s)=
\begin{cases}
a\hat{k}^{\frac{p^{2}-q}{p}}s^{\frac{q-p}{p}}+b, & if s\leq \hat{k},\\
M(s) & if s>k,
\end{cases}
\end{eqnarray*}
is a modified function of $M(s)=as^{p-1}+b$. The positive solutions
of the modified Eq. $(33)$ are critical points of the functional
$$
\mathcal{J}_{\lambda,M_{\hat{k}}}\,(u)=\frac{1}{p}\,\hat{M}_{\hat{k}}\,(\|u\|_{W^{1,p}}^{p})
-\frac{\lambda}{q}\int_{\Omega}f|u|^{q}\,dx-\frac{1}{r}\,\int_{\Omega}g|u|^r\,dx,
$$
where $\hat{M}_{\hat{k}}(s)=\int_{0}^{s}M_{\hat{k}}(t)\,dt$. Note
that by $(4)$, if $u\in \mathcal{N}_{\lambda,M_{\hat{k}}}$ with
$\|u\|_{W^{1,p}}^{p}\leq\hat{k}$, it can be deduced that
\begin{align*}
I''_{u,M_{\hat{k}}}\,(1)&=-a(r-q)\hat{k}^{\frac{p^{2}-q}{p}}\,\|u\|_{W^{1,p}}^{q}\,-b(r-p)\,\|u\|_{W^{1,p}}^{p}+\lambda\,(r-q)\,\int_{\Omega}f|u|^{q}\,dx
\\&\leq-a(r-q)\hat{k}^{\frac{p^{2}-q}{p}}\,\|u\|_{W^{1,p}}^{q}\,-b(r-p)\,\|u\|_{W^{1,p}}^{p}+\lambda\,(r-q)\,S_{q}^{\frac{-q}{p}}\|f^{+}\|_{\infty}\,\|u\|_{W^{1,p}}^{q}
\\&=(r-q)(\lambda\,S_{q}^{\frac{-q}{p}}\|f^{+}\|_{\infty}-a\,\hat{k}^{\frac{p^{2}-q}{p}})\|u\|_{W^{1,p}}^{q}\,-b(r-p)\,\|u\|_{W^{1,p}}^{p}
\end{align*}
\begin{equation}
\hspace*{-6.3 cm}<0,\,\,\text {if}\,\,
\lambda\,\leq\,a\hat{k}^{\frac{p^{2}-q}{p}}\,\|f^{+}\|_{\infty}^{-1}S_{q}^{\frac{q}{p}},
\end{equation}
and for $u\in \mathcal{N}_{\lambda,M_{\hat{k}}}$ with
$\|u\|_{W^{1,p}}^{p}\geq\,\hat{k}$,
\begin{eqnarray*}
I''_{u,M_{\hat{k}}}(1)=[a(p^{2}-r)\,\|u\|_{W^{1,p}}^{p^{2}-r}\,-b(r-q)\,\int_{\Omega}f|u|^{q}\,dx]
\end{eqnarray*}
\begin{equation}
\hspace*{-2 cm}\geq\lambda\,(r-q)\,\int_{\Omega}f|u|^{q}\,dx>0.
\end{equation}
it follows that, by $(33)$,
$\mathcal{N}_{\lambda,M_{\hat{k}}}^{0}\,=\mathcal{N}_{\lambda,M_{\hat{k}}}^{+}\,=\emptyset$
for all
$\lambda\leq\,a\hat{k}^{\frac{p^{2}-q}{p}}\|f^{+}\|_{\infty}^{-1}S_{q}^{\frac{q}{p}}$
and $u\in \mathcal{N}_{\lambda,M_{\hat{k}}}$ with
$\|u\|_{W^{1,p}}^{p}\leq\,\hat{k}$; while by
$(34)$,$\mathcal{N}_{\lambda,M_{\hat{k}}}^{0}\,=\mathcal{N}_{\lambda,M_{\hat{k}}}^{-}\,=\emptyset$
for all $u\in  \mathcal{N}_{\lambda,M_{\hat{k}}}$ with
$\|u\|_{W^{1,p}}^{p}\geq\,\hat{k}$. consequently, if
$0<\lambda\leq\,a\hat{k}^{\frac{p^{2}-q}{p}}\|f^{+}\|_{\infty}^{-1}\,S_{q}^{\frac{q}{p}}$,
the following results are obtained.\\\\
{\bf Lemma 7.1.}
\begin{itemize}
\item[(i)]$\mathcal{N}_{\lambda,M_{\hat{k}}}=
\mathcal{N}_{\lambda,M_{\hat{k}}}^{+}\cup\,
\mathcal{N}_{\lambda,M_{\hat{k}}}^{-}$
$(i.e.\mathcal{N}_{\lambda,M_{\hat{k}}}^{0}=\emptyset)$;
\item[(ii)]$\mathcal{N}_{\lambda,M_{\hat{k}}}\cap\,\{u\in
W_{0}^{1,p}(\Omega):\|u\|_{W^{1,p}}^{p}=\hat{k}\}=\emptyset$;
\item[(iii)]$N_{\lambda,M_{\hat{k}}}^{-}\subset\{u\in
W_{0}^{1,p}(\Omega):\|u\|_{W^{1,p}}^{p}<\hat{k}\}$;
\item[(iv)]$\mathcal{N}_{\lambda,M_{\hat{k}}}^{+}\subset\{u\in
W_{0}^{1,p}(\Omega):\|u\|_{W^{1,p}}^{p}\,>\hat{k}\}$.
\end{itemize}
\hspace{0.6 cm} It is well known that the minimum problem
\begin{equation}
\mathbb{S}=\inf_{u\in M}\,K(u)>0
\end{equation}
can be achieved at positive function $u_{0}\in M$ such that
$K(u_{0})=\mathbb{S}$ ( see \cite {mw}), where
$K(u)=\frac{1}{p}\,\|u\|_{W^{1,p}}^{p}\,-\frac{1}{r}\,\int_{\Omega}g|u|^{r}\,dx$
and $$M=\{u\in
W_{0}^{1,p}(\Omega)\backslash\{0\}:\|u\|_{W^{1,p}}^{p}=\int_{\Omega}g|u|^{r}\,dx\}.$$
Let
$v_{0}=\frac{\hat{k}^{\frac{1}{p}}\,u_{0}}{\|u_{0}\|_{W^{1,p}}}$.
Then $\|v_{0}\|_{W^{1,p}}^{p}=\hat{k}$ and
$$
\int_{\Omega}g|v_{0}|^{r}\,dx=\hat{k}^{\frac{r}{p}}\|u_{0}\|_{W^{1,p}}^{p-r}=\hat{k}^{\frac{r}{p}}(\frac{r-p}{pr\mathbb{S}})^{\frac{(r-p)}{p}}\,>\frac{pb^{2}\,(r-p)}{a(2p-r)^{2}},
$$
provided that
$a<A_{*}=\frac{p^{\frac{r}{(p-r)}}(r-p)^{p}}{r\mathbb{S}}\,(\frac{2p-r}{b})^{\frac{(p^{2}-r)}{(r-p)}}$.\\\\
{\bf Lemma 7.2.} For each $a<A_{*}$ there exists
$0<\tilde{\lambda}_{*}\,\leq\,a\hat{k}^{\frac{p^{2}-q}{p}}\,\|f^{+}\|_{\infty}^{-1}\,S_{q}^{\frac{q}{p}}$
such that for $\lambda<\lambda_{*}$ there exists
$\bar{t_{\lambda}}>1$ such that $\bar{t_{\lambda}}\,v_{0}\,\in
\mathcal{N}_{\lambda,M_{\hat{k}}}^{+}$.\\\\
{\bf Proof.} Let
\begin{eqnarray*}
\bar{m}(\lambda,t)\,&=&a\,t^{2p-r}\|v_{0}\|_{W^{1,p}}^{2p}+bt^{p-r}\,\|v_{0}\|_{W^{1,p}}^{p}-t^{q-r}\,\lambda\,\int_{\Omega}f|v_{0}|^{q}\,dx
\\&=&a\,t^{2p-r}\,\hat{k}^{2}+b\,t^{p-r}\,\hat{k}-t^{q-r}\,\lambda\,\int_{\Omega}f|v_{0}|^{q}\,dx
\,\,\text {for} t>0.
\end{eqnarray*}
Clearly, $\bar{m}(\lambda,t)\to -\infty$ as $t\to 0^{+}$ and
$\bar{m}(\lambda,t)\to \infty$ as $t\to \infty$. Since
$$
\bar{m}'(0,t)=b(r-p)\,\hat{k}\,t^{p-r-1}\,(t^{p}-1),
$$
$\bar{m}'(0,t)=0$ at $t=1$,$\bar{m}',(0,t)<0$ for $t\in (0,1)$ and
$\bar{m}'(0,t)>0$ for $t\in (1,\infty)$. Then $\bar{m}(0,t)$
achieves its minimum at $1$, decreasing for $t\in (0,1)$ and
increasing for $t\in (1,\infty)$. Thus,
\begin{eqnarray*}
\min_{t>0}\,\bar{m}(0,t)=\bar{m}(0,1)=\frac{pb^{2}\,(r-p)}{a(2p-r)^{2}}<\int_{\Omega}g|v_{0}|^{r}\,dx,
\end{eqnarray*}
it follows that there exists $\bar{t_{0}}>1$ such that
\begin{eqnarray*}
\bar{m_{0}}(0,\bar{t_{0}})=\int_{\Omega}g|v_{0}|^{r}\,dx\,\,\text
{and}  \bar{m}_{0}'(0,\bar{t}_{0})>0.
\end{eqnarray*}
By the implicit function theorem, there exists a positive number
$$
\tilde{\lambda}_{*}<a\hat{k}^{\frac{p^{2}-q}{p}}\,\|f^{+}\|_{\infty}^{-1}\,S_{q}^{\frac{q}{p}},
$$
such that for every $\lambda<\tilde{\lambda}_{*}$ there exists
$\bar{t_{\lambda}}>1$ such that
\begin{eqnarray*}
\bar{m}'(\lambda,\bar{t_{\lambda}})=0,
\,\,\,\,\bar{m}(\lambda,\bar{t_{\lambda}})=\int_{\Omega}g|v_{0}|^{r}\,dx.
\end{eqnarray*}
Now,
\begin{align*}
\langle \mathcal{J} _{\lambda
}'(\bar{t_{\lambda}}v_{0}),\bar{t_{\lambda}}v_{0}\rangle
&=a\bar{t}_{\lambda}^{2p}\,\|v_{0}\|_{W^{1,p}}^{2p}+b\,\bar{t}_{\lambda}^{p}\,\|v_{0}\|_{W^{1,p}}^{p}-\bar{t}_{\lambda}^{q}\,\lambda\,\int_{\Omega}f|v_{0}|^{q}\,dx-\bar{t}_{\lambda}^{r}\,\int_{\Omega}g|v_{0}|^{r}\,dx
\\&=\bar{t}_{\lambda}^{r}[\bar{m}(\lambda,\bar{t}_{\lambda})-\int_{\Omega}g|v_{0}|^{r}\,dx]=0
\,\,(i.e.\bar{t}_{\lambda}v_{0}\in
\mathcal{N}_{\lambda,M_{\hat{k}}}),
\end{align*}
and
$$\|\bar{t_{\lambda}}v_{0}\|_{W^{1,p}}^{p}=\bar{t}_{\lambda}^{p}\,\hat{k}>\hat{k}.$$
Thus, by Lemma 7.1\,(iV), $\bar{t}_{\lambda}\,v_{0}\in
\mathcal{N}_{\lambda,M_{\hat{k}}}^{+}$.$\hspace* {0.3cm}\Box$\\\\
{\bf Theorem 7.3.} For each $a<A_{*}$ there exists
$0<\tilde{\lambda}_{*}\,\leq\,a\hat{k}^{\frac{p^{2}-q}{p}}\,\|f^{+}\|_{\infty}^{-1}\,S_{q}^{\frac{q}{p}}$
such that for $0<\lambda<\lambda_{*}$,Eq. $(2)$ has a positive
solution $\hat{u}_{\lambda}$ with
$\|\hat{u}_{\lambda}\|_{W^{1,p}}^{p}>\hat{k}$.\\\\
{\bf Proof.} By Lemma 3.2\,(ii) and Lemma 7.1\,(iV), we have the
energy functional $\mathcal{J}_{\lambda,\hat{k}}$ which is coercive
and bounded below on $\mathcal{N}_{\lambda,M_{\hat{k}}}^{+}$. Then
the minimum problem $\alpha_{\lambda,M_{\hat{k}}}^{+}=\inf_{u\in
\mathcal{N}_{\lambda,M_{\hat{k}}}^{+}}\,\mathcal{J}_{\lambda,M_{\hat{k}}}(u)$
is well defined. Moreover, by $(6)$, Lemma 7.1\,(iV) and
$\mathcal{N}_{\lambda,M_{\hat{k}}}^{+}$ is a nonempty natural
constraint. Thus, by the Ekeland variational principle \cite {ie},
there exists a minimizing sequence $\{u_{n}\}$ for
$\mathcal{J}_{\lambda,M_{\hat{k}}}$ on
$\mathcal{N}_{\lambda,M_{\hat{k}}}^{+}$ such that\\
$$
\mathcal{J}_{\lambda,M_{\hat{k}}}\,(u_{n})=\alpha_{\lambda,{\hat{k}}}^{+}+\,o(1)
\,\,\text {and}
 \mathcal{J}'_{\lambda,M_{\hat{k}}}\,(u_{n})=\,o(1)\,\,\text
 {in}  W^{-1,p'}.
$$
Using a similar argument to that in Lemma $3.1$, there exists
subsequence $\{u_{n}\}$ and ${\hat{u_{\lambda}}\in
W_{0}}^{1,p}\,(\Omega)$ is a nonzero solution of Eq. $(32)$ such
that $u_{n}\to \hat{u}_{\lambda}$ strongly in
$W_{0}^{1,p}\,(\Omega)$, with
$\|\hat{u}_{\lambda}\|_{W^{1,p}}^{p}\,\geq\,\hat{k}$ and
$\mathcal{J}_{\lambda,\hat{k}}\,(\hat{u}_{\lambda})=\alpha_{\lambda,\hat{k}}^{+}$.
Since
$\mathcal{J}_{\lambda,M_{\hat{k}}}\,(\hat{u}_{\lambda})=\mathcal{J}_{\lambda,M_{\hat{k}}}\,(|\hat{u}_{\lambda}|)$,
by Lemma $2.1$, we may assume that $u_{0}$ is a positive solution of
Eq. $(34)$. By lemma $7.1$$(i)$
$\|\hat{u}_{\lambda}\|_{W^{1,p}}^{p}>\hat{k}$. Thus,
$\hat{u}_{\lambda}$ is a positive solution of Eq.
$(2)$.$\hspace* {.2cm}\Box$\\

Now, we complete the proof of Theorem $2.4.$ For each $\theta>0$ and
$0<a<\{\frac{b(r-p)}{rA_{0}L(\theta)},A_{*}\}$, there exists a
positive number
$\tilde{\lambda}_{*}\leq\min\{\theta,\hat{\Lambda}\}$ such that for
$0<\lambda<\tilde{\lambda}_{*}$, Eq. $(2)$ has three positive
solutions $u_{\lambda,M}^{(1),+}$,$u_{\lambda,M}^{(p),+}$ and
$u_{\lambda,M}^{-}$ such that $u_{\lambda,M}^{(i),+}\in
\mathcal{N}_{\lambda,M}^{+}$ and $u_{\lambda,M}^{-}\in
\mathcal{N}_{\lambda,M}^{-}$. Moreover, by $k\in
(\frac{b(r-p)}{ra},\frac{b(r-p)}{pa})$, $(30)$-$(31)$ and Theorem
$2.3$ and $7.3$, we can conclude that
$$\|u_{\lambda,M}^{(1),+}\|_{W^{1,p}}^{p}<(\frac{p\lambda\,(r-q)\,\|f\|_{\infty}s_{q}^{\frac{-q}{p}}}{b(r-p)^{p}})^{\frac{p}{(p-q)}},$$
$$S_{r}^{\frac{r}{p(r-p)}}\,
(\frac{pb(r-1)(p-q)}{r(r-q)\|g^{+}\|_{\infty}})^{\frac{1}{r-p}}<
\|u_{\lambda,M}^{-}\|_{W^{1,p}}^{p}<\frac{b(r-p)}{pa}<\|u_{\lambda,M}^{(p),+}\|_{W^{1,p}}^{p},$$
also this complete the proof of Theorem $2.4$.\\\\

\end{document}